\documentclass[
 aip,
 cha,
 preprint,
 amsmath,
 amssymb,
 floatfix
]{revtex4-2}

\usepackage[utf8]{inputenc}
\usepackage[T1]{fontenc}
\usepackage{graphicx}
\usepackage{amsthm}
\usepackage{multirow}
\usepackage{amsfonts}
\usepackage{mathrsfs}
\usepackage{xcolor}
\usepackage{textcomp}
\usepackage{booktabs}
\usepackage{algorithm}
\usepackage{algorithmicx}
\usepackage{algpseudocode}
\usepackage{listings}
\usepackage[caption=false]{subfig}
\usepackage[colorlinks=true,allcolors=blue]{hyperref}
\allowdisplaybreaks
\newtheorem{thm}{Theorem}[section]
\newtheorem{deff}{Definition}[section]

\numberwithin{deff}{section}
\numberwithin{thm}{section}

\begin{document}

\title{Maps of q-deformed fractional order: From circle to cardioid via crescent}

\author{Sachin Bhalekar}
\email{sachinbhalekar@uohyd.ac.in}
\thanks{Corresponding author}
\affiliation{School of Mathematics and Statistics, University of Hyderabad, Hyderabad, 500046 India}

\author{Prashant M. Gade}
\email{prashant.m.gade@gmail.com}
\affiliation{Ramniranjan Jhunjhunwala College of Arts, Science and Commerce, Ghatkopar West, Mumbai, 400086 India}

\date{\today}

\begin{abstract}
We introduce a class of \(q\)-deformed fractional order maps by replacing the classical binomial memory kernel in discrete fractional dynamics with Gaussian (\(q\)-) binomial coefficients. The proposed framework interpolates between memoryless discrete maps and classical fractional order maps, unifying circle- and cardioid-shaped stability regions through intermediate crescent geometries. Using the \(Z\)-transform and the \(q\)-binomial theorem, we derive characteristic equations and determine the associated stability regions in the complex plane. We further analyze the asymptotic behavior of the memory kernels, showing that the classical fractional kernel exhibits power-law decay, whereas the \(q\)- and \((p,q)\)-deformed kernels exhibit exponential-type localization. The theory is extended to nonlinear logistic-type maps and to a broader \((p,q)\)-deformed framework, where the regime \(p>q\) yields decaying memory kernels and stable dynamics. Numerical simulations illustrate the theoretical results and the interplay between the deformation parameters, memory effects, and stability geometry.
\end{abstract}

\maketitle

\begin{quotation}
\bfseries
Fractional order maps provide a compact way to model discrete-time systems with memory, but their long-range memory kernels often lead to stability regions that differ sharply from those of ordinary maps. This work introduces a deformation of the fractional memory kernel using Gaussian, or \(q\)-, binomial coefficients. The resulting maps continuously connect memoryless dynamics with classical fractional dynamics: their stability regions deform from circles to crescents and then to cardioid-like domains. The additional \((p,q)\)-deformation further enriches the memory structure by allowing exponentially localized kernels. These results provide a geometric and analytic framework for understanding how deformed memory kernels influence stability in discrete nonlinear dynamics.
\end{quotation}
\vspace{-1.0em}

%%%%%%%%%%%%%%%%%%%%%%%%%%%%%%%%%%%%%%%%%%%%%%%%%%%%%%%%%%%%%%%%%%%%%%%%%%%%%%%%%%%%%%%%%%%%%%%%%%

\section{Introduction}

The concept of \(q\)-deformation, introduced by Jackson in the early twentieth century, has become an important tool in several branches of mathematics and physics, including special functions, combinatorics, quantum groups, and dynamical systems \cite{jackson1909q,kac2002quantum,gasper2004basic,koepf1998hypergeometric}. In this work, we extend the ideas of \(q\)-deformation to discrete fractional calculus by introducing a class of \(q\)-deformed fractional order maps.

Jackson defined the \(q\)-derivative of a function \(f(x)\) by
\[
D_q(f(x))
=
\frac{f(qx)-f(x)}{qx-x},
\]
which reduces to the ordinary derivative in the limit \(q\to1\) \cite{ernst2000history}. Consequently,
\[
D_q(x^n)
=
\frac{1-q^n}{1-q}x^{n-1}
=
[n]_q x^{n-1},
\]
where
\[
[n]_q
=
\frac{1-q^n}{1-q}
\]
is known as the \(q\)-analogue of a number $n$. This naturally leads to the \(q\)-factorial
\[
[n]_q!
=
[n]_q[n-1]_q\cdots[1]_q,
\]
and the \(q\)-exponential function
\[
e_q(x)
=
\sum_{n=0}^{\infty}
\frac{x^n}{[n]_q!},
\]
which is an eigenfunction of the \(q\)-derivative \cite{yamano2002some,gasper2004basic}.

Using the \(q\)-factorial, one can define the Gaussian (or \(q\)-) binomial coefficient \cite{konvalina2000unified}. Since these coefficients reduce to the classical binomial coefficients in the limit \(q\to1\), they provide a natural mechanism for introducing deformation into discrete fractional systems while retaining the classical theory as a limiting case.

Fractional calculus extends differentiation and summation to non-integer orders and has emerged as a powerful framework for modelling systems with memory and hereditary effects \cite{samko1993integration,podlubny1998fractional,kilbas2006theory,hilfer2000applications,mainardi2022fractional}. In discrete fractional dynamics, memory is incorporated through weighted contributions from past states, with foundational developments in fractional differences and transform methods given in \cite{miller1989fractional,atici2007transform,atici2009initial,goodrich2015discrete}. Fractional order maps and related coupled-map systems have also been investigated by Bhalekar and co-workers in connection with stability, synchronization, control and nonlinear dynamics \cite{gade2021fractional,pakhare2022synchronization,bhalekar2022stability,joshi2023controlling}. The fractional order maps are extended even to complex order recently\cite{bhalekar2022stabilityc} and spatially extended version of coupled fractional order maps is investigated\cite{bhalekar2022stability}. Classical fractional order maps are characterized by binomial memory kernels whose asymptotic behavior follows a power law, producing long-range memory effects. A natural question is how the memory structure and the associated stability properties change when the underlying binomial coefficients are replaced by their \(q\)-analogues.

In this work, we answer this question by introducing \(q\)-deformation directly into the fractional memory kernel through Gaussian binomial coefficients. This approach differs fundamentally from earlier studies on fractional versions of \(q\)-deformed maps, such as H{\'e}non map, Stefanski map and logistic map are studied, where one first constructs a \(q\)-deformed map and subsequently introduces fractional order\cite{malik2025dynamical,ran2022dynamics,luo2021fractional}. In this context, the deformation of a function and its fractional variant is also considered\cite{bhalekar2025dynamical}. In this work, the deformation acts on the memory kernel itself, leading to a new class of \(q\)-fractional maps.

Fractional \(q\)-calculus and related \(q\)-fractional equations provide another route to incorporating deformation and memory \cite{annaby2012q}. Fractional \(q\)-deformed chaotic maps based on weight-function approaches have previously been investigated in \cite{wu2020fractional}. In that work, the emphasis was primarily on the construction of \(q\)-fractional chaotic maps and their numerical dynamical behavior. In contrast, the present work develops an analytic framework for \(q\)- and \((p,q)\)-deformed fractional maps using \(Z\)-transform techniques, explicit characteristic equations, and stability geometry in the complex plane.

A central theme of the present work is the interplay between memory kernels and stability geometry. We show that replacing classical binomial coefficients by Gaussian binomial coefficients produces a family of maps whose stability properties depend jointly on the fractional order and the deformation parameter \(q\). The corresponding stability regions evolve continuously from circles associated with classical discrete maps to crescent-shaped domains and finally to the cardioid-shaped regions characteristic of fractional order maps \cite{gade2021fractional,bhalekar2022stability,edelman2023stability}. Thus the \(q\)-deformation provides a geometric bridge between memoryless and fractional dynamics.

The deformation also modifies the asymptotic behavior of the memory kernel itself. While the classical fractional kernel exhibits algebraic power-law decay, the \(q\)-deformed kernel approaches a finite asymptotic value. We further extend the framework to \((p,q)\)-deformations, where an additional deformation parameter allows exponentially localized memory kernels when \(p>1\). Consequently, the proposed framework provides a hierarchy of memory structures ranging from power-law memory to asymptotically saturated memory and exponential memory.

Beyond the single-parameter \(q\)-deformation, the present framework admits a natural extension to \((p,q)\)-deformations. By replacing Gaussian binomial coefficients with \((p,q)\)-binomial coefficients, we obtain a broader class of fractional order maps with richer memory kernels and stability geometries. The \(q\)-deformed maps arise as the special case \(p=1\).

The main contributions of this work are threefold. First, we derive characteristic equations and stability boundaries for \(q\)- and \((p,q)\)-deformed fractional maps using \(Z\)-transform techniques and \(q\)-binomial identities. Second, we demonstrate a continuous geometric transition of the stability region from circles to crescents and finally to cardioids as the deformation parameter varies. Third, we show that the associated memory kernels interpolate between power-law, asymptotically saturated, and exponentially localized memory profiles. We further illustrate the theory through nonlinear logistic-type maps and numerical simulations that support the analytical stability results.

\section{Preliminaries} \label{prel}
In this section, we present some basic definitions and results.\\ 
Let $h>0$, \,$ a\in \mathbb{R}$,
$(h\mathbb{N})_a=\{a,a+h,a+2h,\dots\}$ and $\mathbb{N}_a=\{a,a+1,a+2,\dots\}$.

%\begin{deff}(See\cite{bastos2011discrete,ferreira2011fractional,mozyrska2015transform})
%	For a function $x:(h\mathbb{N})_a\longrightarrow \mathbb{R}$, the forward h-difference operator is defined as
%	$$(\Delta_hx)(t)=\frac{x(t+h)-x(t)}{h},$$
%	where $t\in(h\mathbb{N})_a$. \\
%	Throughout this paper, we take $a=0$ and $h=1$. 
%\end{deff}

%\begin{deff}\cite{mozyrska2015transform}
%	For a function $x:\mathbb{N}_0\longrightarrow \mathbb{R}$, the fractional sum of order $\beta >0$\, is given by 
%	$$(\Delta^{-\beta} x)(t)=\frac{1}{\Gamma(\beta)} \sum_{s=0}^n{\frac{\Gamma(\beta+n-s)}{\Gamma(n-s+1)}}x(s),$$  
%	where $t=\beta +n$, $n\in \mathbb{N}_0$.
%\end{deff}

%\begin{deff}\cite{fulai2011existence,mozyrska2015transform}  
%	Let $\mu>0$ and $m-1<\mu<m$, where $m\in\mathbb{N}.$
%	The $\mu$th fractional Caputo-like difference is defined as
%	$$\Delta^\mu x(t)=\Delta^{-(m-\mu)}(\Delta^mx(t)),$$
%	where $t\in\mathbb{N}_{m-\mu}$ and
%	$$\Delta^mx(t)=\sum_{k=0}^m \left(\begin{array}{c}m\\k\end{array}\right)(-1)^{m-k}x(t+k).$$
%\end{deff}

\begin{deff}\cite{mozyrska2015transform} 
	The Z-transform of a sequence $\{y(n)\}_{n=0}^\infty$ is a complex function given by $Y(z)=Z[y](z)=\sum_{k=0}^\infty y(k)z^{-k}$, where $z\in \mathbb{C}$ is a complex number for which the series converges absolutely.
\end{deff}
\begin{deff}\cite{mozyrska2015transform} 
	Let 
    \begin{eqnarray*}
		\Tilde{\phi}_\alpha(n) &=& \frac{\Gamma(n+\alpha)}{\Gamma(\alpha)\Gamma(n+1)} \nonumber \\
		&=& \left(\begin{array}{c} n+\alpha-1\\n \end{array} \right) = (-1)^n \left( \begin{array}{c}
			-\alpha \\ n \end{array} \right)
	\end{eqnarray*}
    be a family of binomial functions defined on $\mathbb{Z}$, parametrized by $\alpha$.\\
Then, 
	\begin{eqnarray*}
		Z(\Tilde{\phi}_\alpha(t))=\frac{1}{(1-z^{-1})^\alpha},\; |z|>1. \nonumber
	\end{eqnarray*}
\end{deff}

\begin{deff}\cite{mozyrska2015transform}
	The convolution $\phi\ast x$ of the functions $\phi$ and $x$ defined on $\mathbb{N}_0$ is defined as 
	\begin{eqnarray*}(\phi\ast x)(n)=\sum_{s=0}^n \phi(n-s)x(s)=\sum_{s=0}^n \phi(s) x(n-s).
	\end{eqnarray*} 
	Then, the Z-transform of this convolution is
	\begin{eqnarray*}
		Z(\phi\ast x)(n)=(Z(\phi)(n)) Z((x)(n)).
	\end{eqnarray*} 
\end{deff}

\begin{deff} \cite{gasper2004basic}
	Let $a$ and $q$	be real numbers with $|q|<1$ and $n\in \mathbb{N}$. The $q-$Pochhammer symbol  or  $q-$shifted factorial is defined as
	\begin{equation}
		\left(a;q\right)_n=\prod_{k=0}^{n-1}\left(1-aq^k\right).
	\end{equation}
	%This can be defined for the extreme cases $n=0$ and $n=\infty$.
	We have, $\left(a;q\right)_0=1$. Furthermore, this symbol can be extended to an infinite product
	\begin{equation}
		\left(a;q\right)_\infty=\prod_{k=0}^{\infty}\left(1-aq^k\right).
	\end{equation}
	This is analytic function of $q$ for $|q|<1$. 
\end{deff}

\begin{thm}\cite{andrews1986q} ($q-$Binomial Theorem)
For $|z|<1$ and $|q|<1$,
\begin{equation}
\sum_{n=0}^{\infty}\frac{\left(a;q\right)_n}{\left(q;q\right)_n}z^n=\frac{\left(az;q\right)_\infty}{\left(z;q\right)_\infty}. \label{qbin}
\end{equation}
\end{thm}

\begin{deff}\cite{gasper2004basic,weisstein2017binomial,mukhinsymmetric}
The Gaussian binomial coefficient is defined as
\begin{equation}
\binom{m}{r}_q =\prod_{k=1}^{r} \frac{1-q^{m-r+k}}{1-q^{k}}, \label{gb}
\end{equation}
where $m$, $r$ $\in \mathbb{N}_0$.
\end{deff}
Note: We may allow $m$ to take positive real values in (\ref{gb}).
%\begin{thm}\cite{fulai2011existence}
%	The difference equation $$\Delta^{\alpha} x(t)=f(x(t+\alpha-1))-x(t+\alpha-1),$$ where $0<\alpha \leq 1$, $t \in \mathbb{N}_{1-\alpha}$, 
%	is  equivalent to
%	\begin{equation} 
%		x(t)=x(0)+\frac{1}{\Gamma(\alpha)} \sum_{j=1}^t \frac{\Gamma(t-j+\alpha)}{\Gamma(t-j+1 )} (f(x(j-1))-x(j-1)), \label{1}
%	\end{equation}
%	where $t \in \mathbb{N}_{0}$.
%\end{thm}

%\begin{deff}\cite{elaydi2005systems,hirsch2012differential}
%	A steady state solution or an equilibrium $x_*$ of (\ref{1}) is a real number satisfying $f(x_*) = x_*$.
%\end{deff}

%\begin{deff} \cite{elaydi2005systems,hirsch2012differential}
%	An equilibrium $x_*$ is stable if for each  $\epsilon>0$, there exists $\delta>0$ such that $|x_0 -x_* | < \delta $ implies
%	$|x(t) - x_* | < \epsilon$, $t=1,2,3,...$\\
%	If $x_*$ is not stable then it is unstable.
%\end{deff}

%\begin{deff} \cite{elaydi2005systems,hirsch2012differential}
%	An equilibrium point $x_*$  is asymptotically stable  if it is stable and there exists $\delta>0$ such that $|x_0 -x_* | < \delta $ implies $ lim_{t\to\infty}x(t) =x_*$.
%\end{deff}
\subsection{Fractional order maps}
The stability analysis of the following fractional order map 
\begin{eqnarray}
	x(t+1)&=&x(0)+\sum_{j=0}^t \binom{t-j+\alpha -1}{t-j} (a-1)x(j),\nonumber\\
	&=&x(0)+ (a-1) \left(\Tilde{\phi}_\alpha\ast x\right)(t)  \label{2}
\end{eqnarray}
where $0<\alpha<1$, $a$ is a constant and $x(t)\in \mathbb{R} $ for each $t=0,1,2,\cdots$ is discussed in \cite{gade2021fractional,bhalekar2022stability,bhalekar2023fractional}. 
\par In this work, we generalize the map (\ref{2}) by replacing $\Tilde{\phi}_\alpha(t)$ with the following $q-$variant.
\begin{equation}
\Tilde{\phi}_{\alpha,q}(t)=\binom{t+\alpha-1}{t}_q, \label{qb1}
\end{equation}
where $|q|<1$, $0<\alpha<1$, and $t\in \mathbb{N}_{0}$.\\
We can write the generalized $q-$binomial coefficient $\Tilde{\phi}_{\alpha,q}(t)$  in terms of $q-$Pochhammer symbol as below.
\begin{eqnarray}
\Tilde{\phi}_{\alpha,q}(t) &=& \prod_{k=1}^{t} \frac{1-q^{\alpha+k-1}}{1-q^{k}}
= \frac{\prod_{k=0}^{t-1}\left(1-q^{\alpha}q^k\right)}{\prod_{k=0}^{t-1}\left(1-q q^{k}\right)}\nonumber\\
&=& \frac{\left(q^{\alpha};q\right)_t}{\left(q;q\right)_t}.\label{qb2}
\end{eqnarray}
We can find the Z-transform of this function by using the q-binomial theorem (\ref{qbin}) as below:
\begin{eqnarray}
Z\left(\Tilde{\phi}_{\alpha,q}(t)\right) &=&\sum_{n=0}^{\infty}\frac{\left(q^{\alpha};q\right)_n}{\left(q;q\right)_n}z^{-n}\nonumber\\
&=&\frac{\left(q^{\alpha}z^{-1};q\right)_\infty}{\left(z^{-1};q\right)_\infty}. \label{zt1}
\end{eqnarray}
\section{\texorpdfstring{The $q$-variant of fractional order map}{The q-variant of fractional order map}}
We consider the map
\begin{eqnarray}
	x(t+1)&=&x(0)+\sum_{j=0}^t \binom{t-j+\alpha -1}{t-j}_q (a-1)x(j),\nonumber\\
	&=&x(0)+ (a-1) \left(\Tilde{\phi}_{\alpha,q}\ast x\right)(t),  \label{3}
\end{eqnarray}
where $|q|<1$, $0<\alpha<1$, and $t\in \mathbb{N}_{0}$.\\
In order to discuss the stability of this map, we need to find the characteristic equation. Taking the $Z$-transform of (\ref{3}) and using (\ref{zt1}) and elementary properties, we get
\begin{equation}
z X(z)-z x(0)=\frac{x(0)}{1-z^{-1}}+(a-1)\frac{\left(q^{\alpha}z^{-1};q\right)_\infty}{\left(z^{-1};q\right)_\infty}X(z),
\label{eq10}
\end{equation}
where $X(z)$ denotes the  $Z$-transform of the sequence $x(t)$. Thus, the characteristic equation of our system  (\ref{3}) is
\begin{equation}
a=1+\frac{z \left(z^{-1};q\right)_\infty}{\left(q^{\alpha}z^{-1};q\right)_\infty}. \label{cheq1}
\end{equation}
The zero solution of the system (\ref{3}) is asymptotically stable if all the roots $z$ of the characteristic equation (\ref{cheq1}) satisfy $|z|<1$. However, checking this condition is difficult. Therefore, we will find the boundary of the stable region, defined by $\left\{z=e^{\iota \theta} | 0\leq\theta<2\pi\right\}$ in terms of the parameter $a$. 
\par If we assume that $a$ is real number, then the boundary points are $z=-1$ and $z=1$. Substituting these values in (\ref{cheq1}), we get the bounds as $a=1-\frac{ \left(-1;q\right)_\infty}{\left(-q^{\alpha};q\right)_\infty}$ and $a=1$, respectively.
Let $\Lambda_{q,\alpha}= \frac{ \left(-1;q\right)_\infty}{\left(-q^{\alpha};q\right)_\infty}$ and we  have following result

\begin{thm}
The zero solution of the fractional order map (\ref{3}) involving $q-$binomial coefficient is locally asymptotically stable if $1-
%\frac{ \left(-1;q\right)_\infty}{\left(-q^{\alpha};q\right)_\infty}
\Lambda_{q,\alpha}
<a<1$. If we allow the parameter $a$ to take complex values then the stability region in the complex plane is bounded by the parametric curve $\gamma(\theta)=1+\frac{e^{\iota \theta} \left(e^{-\iota \theta};q\right)_\infty}{\left(q^{\alpha}e^{-\iota \theta};q\right)_\infty}$, $0\leq\theta<2\pi$. The stability is assured if the parameter $a$ lies inside the positively oriented, simple, bounded curve $\gamma(\theta)$.
\end{thm}

\begin{figure}
	\subfloat[$q=0$]{\includegraphics[width = 2in]{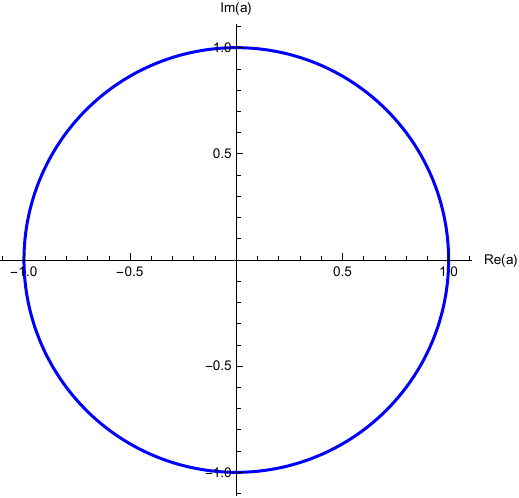}} 
	\subfloat[$q=0.25$]{\includegraphics[width = 2in]{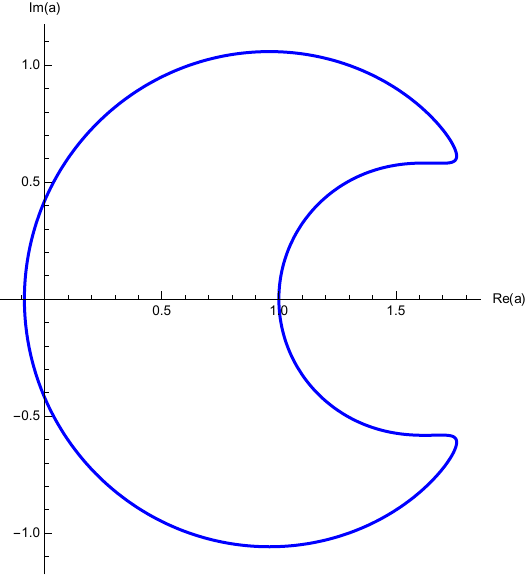}}
	\subfloat[$q=0.4$]{\includegraphics[width = 2in]{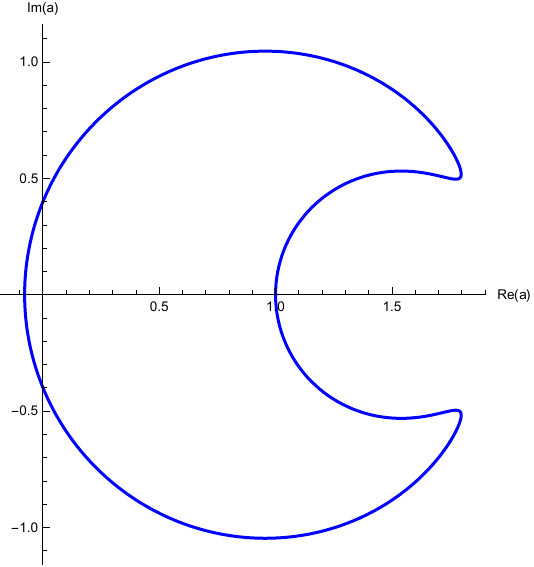}}\\
	\subfloat[$q=0.65$]{\includegraphics[width = 2in]{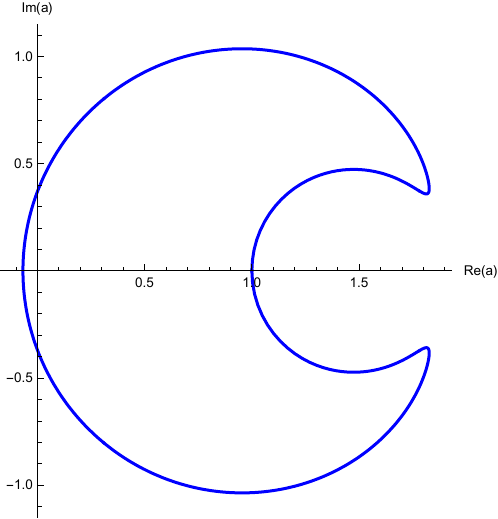}} 
	\subfloat[$q=0.85$]{\includegraphics[width = 2in]{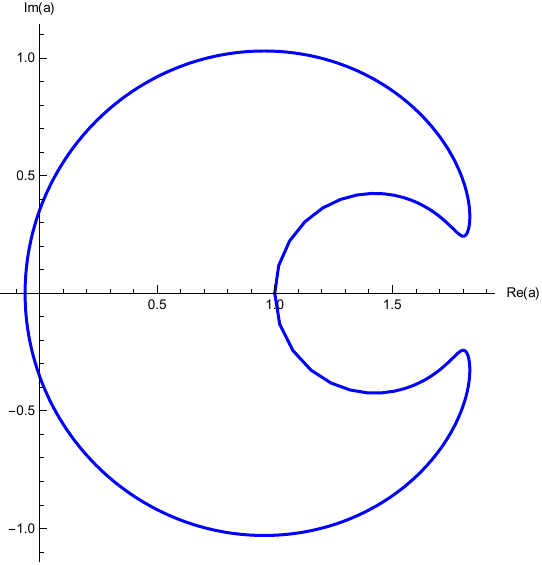}}
	\subfloat[$q=1$]{\includegraphics[width = 2in]{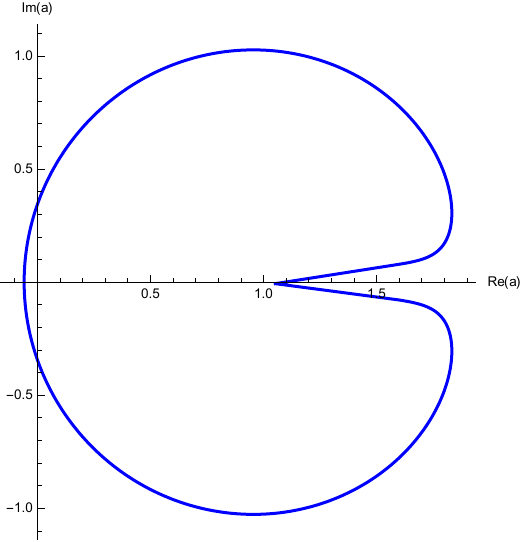}}
	\caption{Stable regions of the system (\ref{3}) for $\alpha=0.082$ and various values of $q$}
	\label{fig0}
\end{figure}

\begin{figure}
	\subfloat[$q=0$]{\includegraphics[width = 2in]{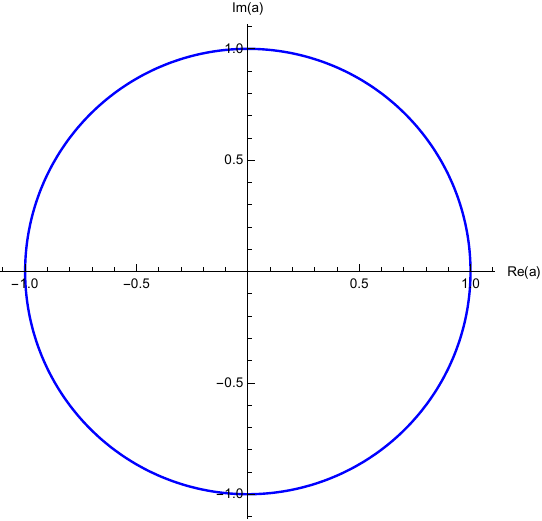}} 
	\subfloat[$q=0.2$]{\includegraphics[width = 2in]{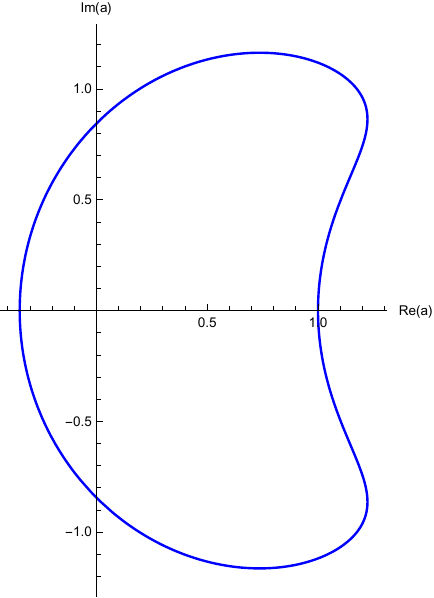}}
	\subfloat[$q=0.4$]{\includegraphics[width = 2in]{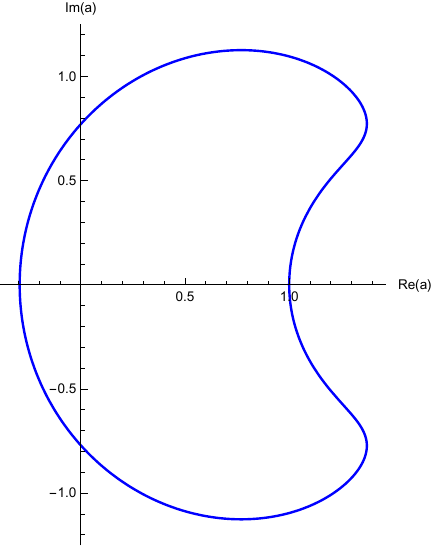}}\\
	\subfloat[$q=0.6$]{\includegraphics[width = 2in]{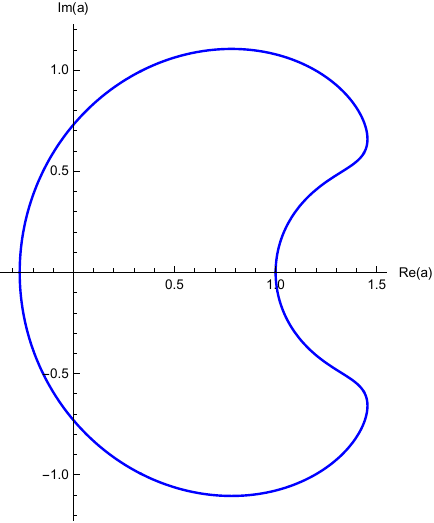}} 
	\subfloat[$q=0.8$]{\includegraphics[width = 2in]{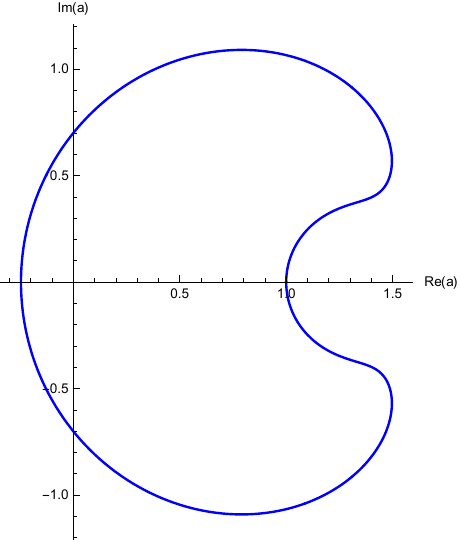}}
	\subfloat[$q=1$]{\includegraphics[width = 2in]{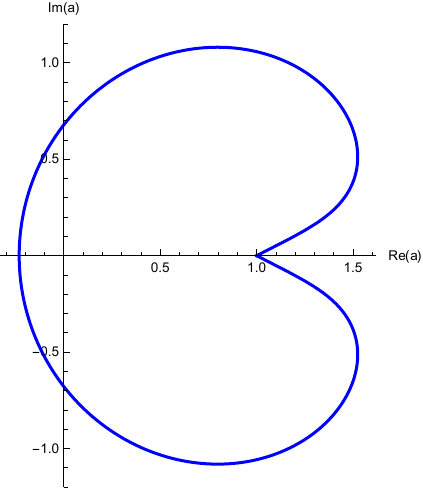}}
	\caption{Stable regions of the system (\ref{3}) for $\alpha=0.3$ and various values of $q$}
	\label{fig1}
\end{figure}

\begin{figure}
	\subfloat[$q=0$]{\includegraphics[width = 2in]{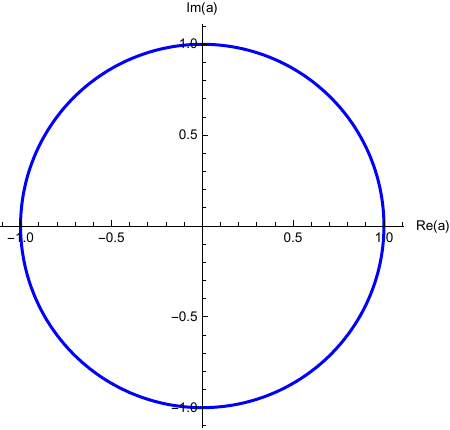}} 
	\subfloat[$q=0.2$]{\includegraphics[width = 2in]{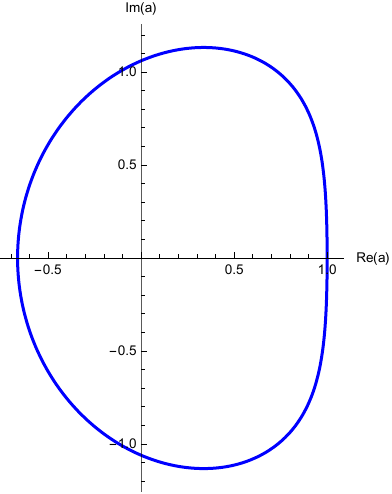}}
	\subfloat[$q=0.4$]{\includegraphics[width = 2in]{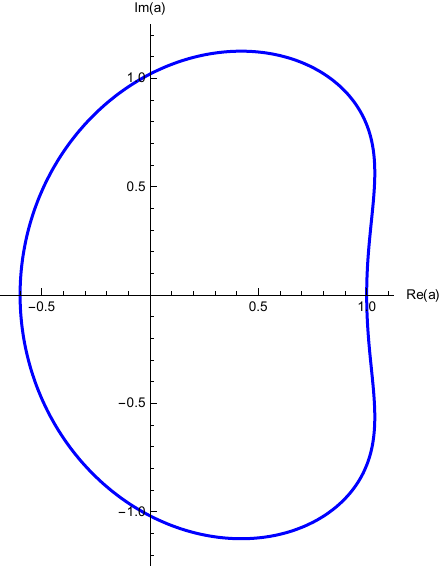}}\\
	\subfloat[$q=0.6$]{\includegraphics[width = 2in]{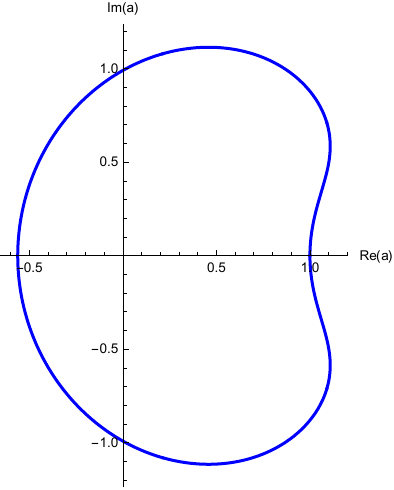}} 
	\subfloat[$q=0.8$]{\includegraphics[width = 2in]{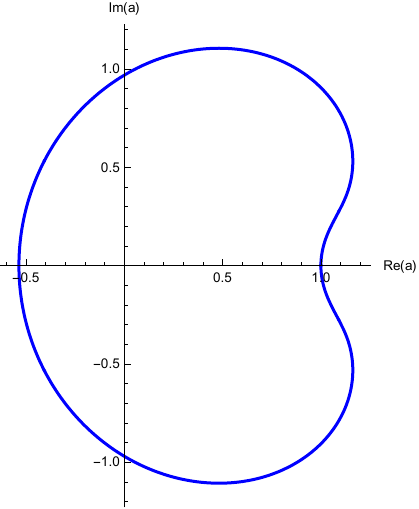}}
	\subfloat[$q=1$]{\includegraphics[width = 2in]{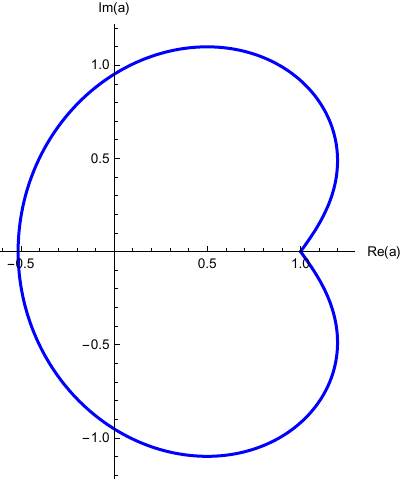}}
	\caption{Stable regions of the system (\ref{3}) for $\alpha=0.6$ and various values of $q$}
	\label{fig2}
\end{figure}

\begin{figure}
	\subfloat[$q=0$]{\includegraphics[width = 2in]{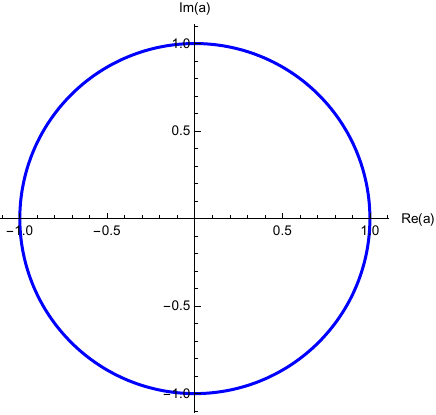}} 
	\subfloat[$q=0.2$]{\includegraphics[width = 2in]{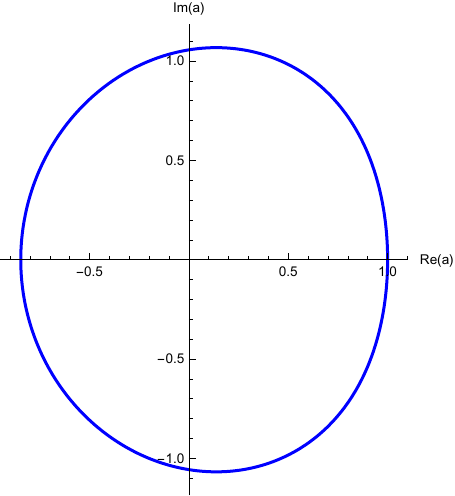}}
	\subfloat[$q=0.4$]{\includegraphics[width = 2in]{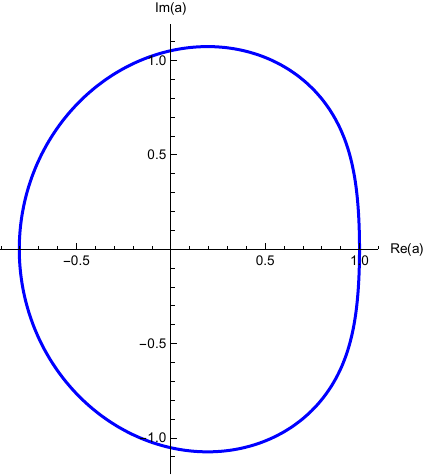}}\\
	\subfloat[$q=0.6$]{\includegraphics[width = 2in]{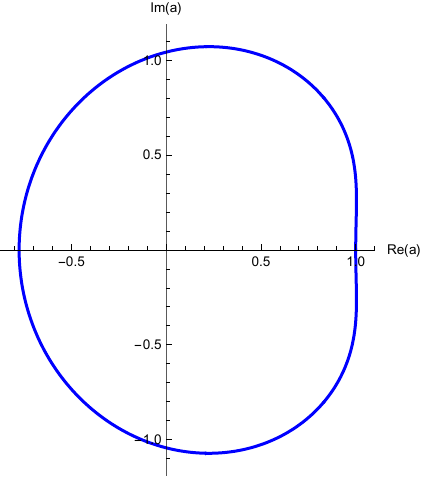}} 
	\subfloat[$q=0.8$]{\includegraphics[width = 2in]{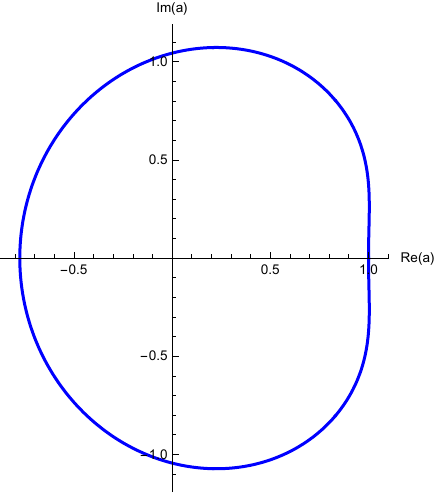}}
	\subfloat[$q=1$]{\includegraphics[width = 2in]{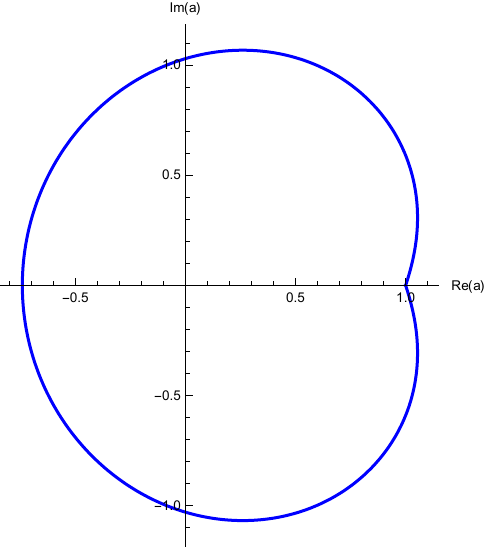}}
	\caption{Stable regions of the system (\ref{3}) for $\alpha=0.8$ and various values of $q$}
	\label{fig3}
\end{figure}

For $0<\alpha<1$, it is observed that, the curve $\gamma(\theta)$ is a circle for $q=0$ and a cardioid for $q=1$. The curve $\gamma(\theta)$ is a crescent for $0<q<1$. The extreme cases,  $q=0$ and  $q=1$  represents the classical map $x(t+1)=ax(t)$ and fractional order map $x(t+1)=x(0)+ (a-1) \left(\Tilde{\phi}_\alpha\ast x\right)(t)$, respectively. These stable regions are well-studied in the literature \cite{gade2021fractional,bhalekar2022stability,edelman2023stability}. Our work unifies these two maps and there is an interplay between circle, cardioid and crescent. We sketch the boundary curve  $\gamma(\theta)$ for various values of $\alpha$ and $q$ in Figures \ref{fig0}, \ref{fig1}, \ref{fig2} and \ref{fig3}.

\subsection{\texorpdfstring{Limiting case $q\rightarrow 1$ and $q\rightarrow 0$}{Limiting case q -> 1 and q -> 0}}
From equation (\ref{cheq1}), the characteristic equation is
\begin{equation}
z-(a-1)\frac{(q^{\alpha}z^{-1};q){\infty}}{(z^{-1};q){\infty}}=0.
\end{equation}

Therefore, by (\ref{zt1})
\begin{equation}
z-(a-1)\sum_{n=0}^{\infty}\frac{(q^{\alpha};q)_n}{(q;q)_n}z^{-n}=0. \label{n1}
\end{equation}

Now, \cite{gasper2004basic}
\begin{equation}
\lim_{q\to 1}\frac{(q^{\alpha};q)_n}{(q;q)_n}
=
\binom{n+\alpha-1}{n}.
\end{equation}

Taking limit of (\ref{n1}) as $q\to 1$, we get
\begin{equation}
z-(a-1)\sum_{n=0}^{\infty}
\binom{n+\alpha-1}{n}z^{-n}=0.
\end{equation}

The binomial theorem thus gives
\begin{equation}
z-(a-1)(1-z^{-1})^{-\alpha}=0.
\end{equation}

Therefore, 
\begin{equation}
a=1+z(1-z^{-1})^{\alpha}.
\end{equation}

This is the characteristic equation of fractional order map (\ref{2}), and it gives
a cardioid \cite{gade2021fractional} for
\begin{equation}
z=e^{i\theta},
\qquad
0\leq \theta \leq 2\pi.
\end{equation}

Thus, $q\to 1$ gives fractional-order map.
\begin{equation}
x(t+1)=x(0)+\sum_{j=0}^{t}
\binom{t-j+\alpha-1}{t-j}(a-1)x(j).
\end{equation}

 On the other hand,
using the definition of the infinite $q$-shifted factorial,
\[
(a;q)_\infty=\prod_{k=0}^{\infty}(1-aq^k),
\]
and for $q\rightarrow 0$

\begin{equation}
\lim_{q\to 0}(q^{\alpha}z^{-1};q)_{\infty}=1.
\label{20}
\end{equation}

Substituting $a=z^{-1}$ gives
\[
(z^{-1};q)_\infty
=
\prod_{k=0}^{\infty}(1-z^{-1}q^k).
\]
Separating the $k=0$ term, we obtain
\[
(z^{-1};q)_\infty
=
(1-z^{-1}q^0)\prod_{k=1}^{\infty}(1-z^{-1}q^k)
=
(1-z^{-1})\prod_{k=1}^{\infty}(1-z^{-1}q^k).
\]
Thus
\begin{equation}
\lim_{q\to 0}(z^{-1};q)_{\infty}
= (1-z^{-1})=
\frac{z-1}{z}.
\label{21}
\end{equation}

Therefore, taking limit of characteristic equation as $q\to 0$, we get
\begin{equation}
z-(a-1)\frac{z}{z-1}=0.
\label{22}
\end{equation}

i.e.
\begin{equation}
a=z=e^{i\theta},
\qquad
0\leq \theta \leq 2\pi.
\end{equation}

This gives a circle centered at origin and radius $1$.

This is stable region for classical case
\begin{equation}
x_{n+1}=ax_n.
\end{equation}

%\begin{figure}
%	\centering\subfloat[$q=0$]{%
%		\includegraphics[scale=0.24]{al3q0}
%	}
%	\hspace{0.5cm}
%	\centering\subfloat[[$q=0.2$]{%
%		\includegraphics[scale=0.2]{al3q2}
%	}
%	\hspace{0.5cm}
%	\centering\subfloat[$q=0.4$]{%
%		\includegraphics[scale=0.22]{al3q4}
%	}
%	\hspace{0.5cm}
%	\centering\subfloat[$q=0.6$]{%
%		\includegraphics[scale=0.21]{al3q6}
%	}\\
%	\centering\subfloat[$q=0.8$]{%
%		\includegraphics[scale=0.21]{al3q8}
%	}
%	\hspace{0.5cm}
%	\centering\subfloat[$q=1$]{%
%		\includegraphics[scale=0.26]{al3q10}
%	}
%	\caption{Stable regions of the system (\ref{3}) for $\alpha=0.3$ and various values of $q$}
%	\label{fig1}
%\end{figure}

\begin{figure}
\includegraphics[width=5in]{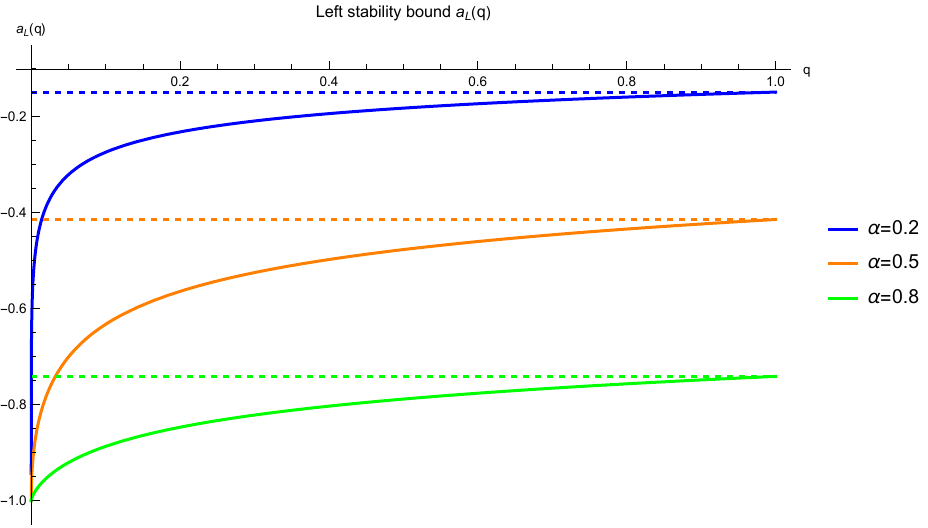}
\caption{
Variation of the left stability bound
$
a_L(q)=1-\frac{(-1;q)_\infty}{(-q^\alpha;q)_\infty}
$
with the deformation parameter $q$ for $\alpha=0.2$, $0.5$, and $0.8$. The bound evolves continuously from $a_L=-1$ at $q=0$ to the fractional-order limit $a_L=1-2^\alpha$ as $q\to1^-$. Increasing $\alpha$ modifies the real stability interval.
}
\end{figure}

\section{Illustrative Examples}
We validated the stable regions presented in the previous section by taking particular values of parameters and solving the system (\ref{3}), following the numerical-validation approach used for fractional maps in \cite{gade2021fractional,pakhare2022synchronization,bhalekar2023fractional}. The results are summarized in Table \ref{tab1}.\\
\begin{table}[h]
\begin{tabular}{|c|c|c|c|c|c|}
	\hline
 $\alpha$ & $q$ & Stable region &  $a$ & Stability property & Solution  \\
	\hline
 $0.2$ &$0.3$  & $(-0.20894, 1)$& $-0.1$ & Stable & Fig. \ref{fig4} (a) \\
	\cline{4-6}
	  &  &  & $-0.22$ & Unstable & Fig. \ref{fig4} (b)\\
	  	\cline{4-6}
	  &  &  & $1.1$ & Unstable & Fig. \ref{fig4} (c)\\
	\hline
 $0.5$ &$0.2$  & $(-0.56385, 1)$  & $-0.5$  &Stable & Fig. \ref{fig4} (d)\\
 	\cline{4-6}
 &  &  & $-0.6$ & Unstable & Fig. \ref{fig4} (e)\\
	\hline
 $0.3$ & $0.9$  & As shown in Fig. \ref{fig5}  & $-0.1-0.2\iota$ & Stable & Fig. \ref{fig4} (f)\\
 \cline{4-6}
 &  &  & $0.1+\iota$ & Unstable & Fig. \ref{fig4} (g)\\
  \cline{4-6}
 &  &  & $1.2+0.12\iota$ & Unstable & Fig. \ref{fig4} (h)\\
	\hline
%	\caption{Validation of stability result}
\end{tabular}
\caption{Validation of stability result}
\end{table}\label{tab1}

\begin{figure}
	\subfloat[$\alpha=0.2, q=0.3, a=-0.1$]{\includegraphics[width = 2in]{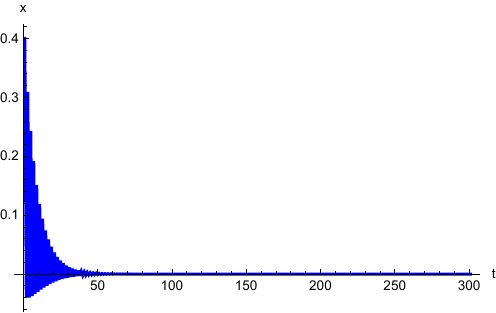}} 
	\subfloat[$\alpha=0.2, q=0.3, a=-0.22$]{\includegraphics[width = 2in]{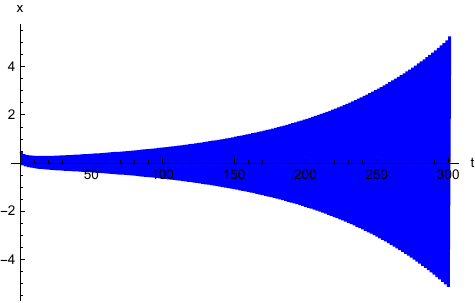}}\\
	\subfloat[$\alpha=0.2, q=0.3, a=1.1$]{\includegraphics[width = 2in]{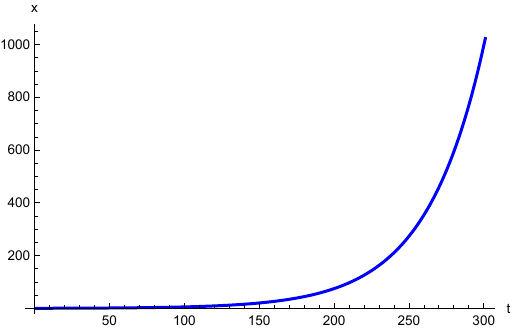}}
	\subfloat[$\alpha=0.5, q=0.2, a=-0.5$]{\includegraphics[width = 2in]{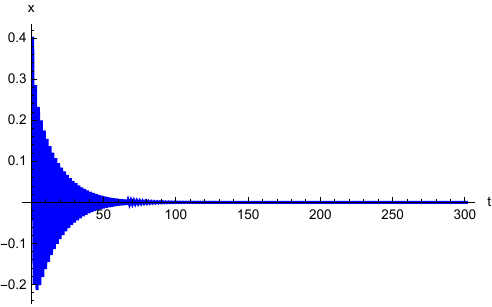}} \\
	\subfloat[$\alpha=0.5, q=0.2, a=-0.6$]{\includegraphics[width = 2in]{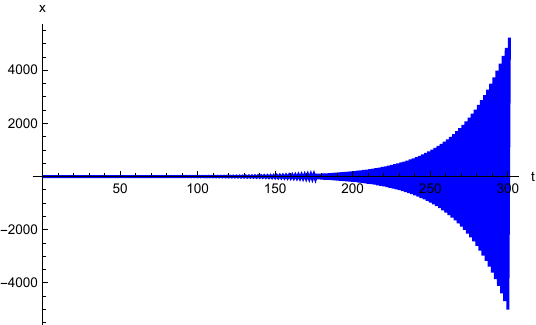}}
	\subfloat[$\alpha=0.3, q=0.9, a=-0.1-0.2\iota$]{\includegraphics[width = 2in]{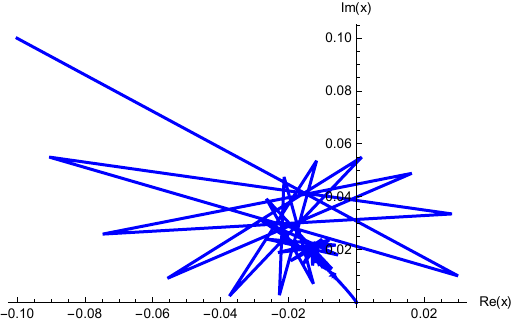}}\\
	\subfloat[$\alpha=0.3, q=0.9, a=0.1+\iota$]{\includegraphics[width = 2in]{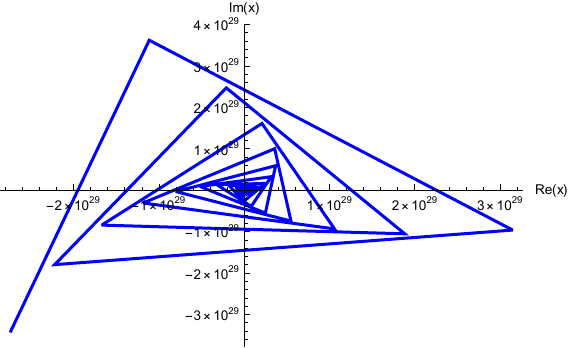}}
		\subfloat[$\alpha=0.3, q=0.9, a=1.2+0.12\iota$]{\includegraphics[width = 2in]{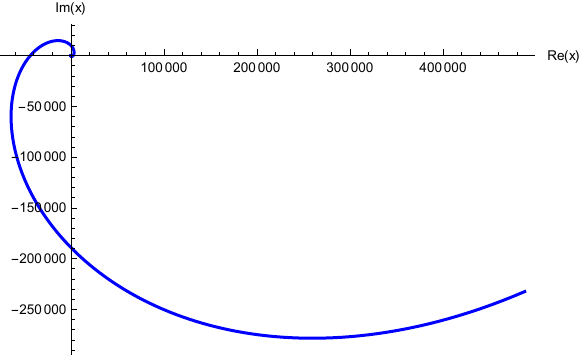}}	
	\caption{Solutions of the system (\ref{3}) for various values of $\alpha$ $q$ and $a$}
	\label{fig4}
\end{figure}

\begin{figure}
\includegraphics{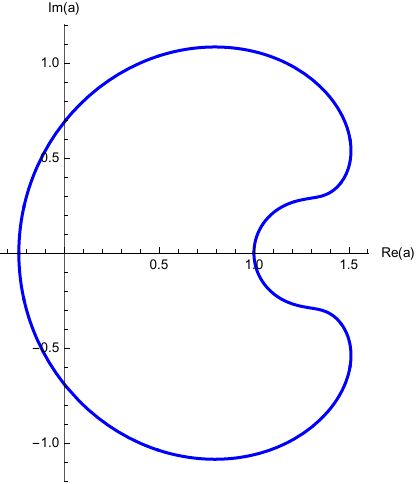}
\caption{Stable region for the system (\ref{3}) for $\alpha=0.3$ and $q=0.9$}
\label{fig5}
\end{figure}

\section{Nonlinear Maps}
Consider the nonlinear version of the system  (\ref{3}) defined by
\begin{equation}
x(t+1)=x(0)+\sum_{j=0}^t \binom{t-j+\alpha -1}{t-j}_q \left(f(x(j))-x(j)\right), \label{nlin}
\end{equation}
where $0<\alpha<1$ and $0<q<1$.
\par The fixed point $x_*$ of the map $f$ is called an equilibrium point of the system (\ref{nlin}). The linearization of the system (\ref{nlin}) near $x_*$ is given by the equation (\ref{3}) with $a=f'(x_*)$.
%Let $\Lambda_{q,\alpha}=\frac{(-1;q)_\infty}{(-q^\alpha;q)_\infty}.$  
The equilibrium  $x_*$ is locally asymptotically stable if 
$1-\Lambda_{q,\alpha}<f'(x_*)<1$.
%$1-\frac{ \left(-1;q\right)_\infty}{\left(-q^{\alpha};q\right)_\infty}
\par We discuss the stablility properties of the equilibrium points of the system (\ref{nlin}) with the logistic map $f(x)=rx(1-x)$, a standard nonlinear test map also considered in fractional-order map studies \cite{wu2014chaos,gade2021fractional,bhalekar2023fractional}. The fixed points of the logistic map are $x_{1*}=0$ and $x_{2*}=1-1/r$. 

Now, $f'(x_{1*})=r$ and $f'(x_{2*})=2-r$. 
%Let $$\Lambda_{q,\alpha}=\frac{(-1;q)_\infty}{(-q^\alpha;q)_\infty}.$$ 
Therefore, the equilibrium $x_{1*}$ is locally asymptotically stable (LAS) for 
$r \in (1-\Lambda_{q,\alpha}, 1)$, whereas the stable region for $x_{2*}$ is $r \in (1, 1+\Lambda_{q,\alpha})$.

%Now, $f'(x_{1*})=r$ and $f'(x_{2*})=2-r$. Therefore, the stable regions %for $x_{1*}$ and $x_{2*}$ are $1-\frac{ \left(-1;q\right)_\infty}{\left(-%q^{\alpha};q\right)_\infty}<r<1$ and $1<r<1+\frac{ %\left(-1;q\right)_\infty}{\left(-q^{\alpha};q\right)_\infty}$, %respectively.

%For the fixed values $\alpha=0.6$ and $q=0.3$, we have $1-\frac{ %\left(-1;q\right)_\infty}{\left(-q^{\alpha};q\right)_\infty}=-0.62685$ %and $1+\frac{ \left(-1;q\right)_\infty}{\left(-%q^{\alpha};q\right)_\infty}=2.62685$. 

For the fixed values $\alpha=0.6$ and $q=0.3$, we have
$\Lambda_{q,\alpha}=1.62685$.
%$1-\Lambda_{q,\alpha}=-0.62685$ and $1+\Lambda_{q,\alpha}=2.62685$.
Therefore, the equilibrium $x_{1*}$ is locally asymptotically stable (LAS) for $r\in (-0.62685, 1)$ whereas the stable region for $x_{2*}$ is $r\in (1, 2.62685)$.

 %Let $$\Lambda=\frac{(-1;q)_\infty}{(-q^\alpha;q)_\infty}.$$ Therefore, %the equilibrium $x_{1*}$ is locally asymptotically stable (LAS) for 
%$r \in (1-\Lambda, 1)$, whereas the stable region for $x_{2*}$ is $r \in 
%(1, 1+\Lambda)$.

\begin{figure}
	\subfloat[$r=0.9, x(0)=0.8$]{\includegraphics[width = 2in]{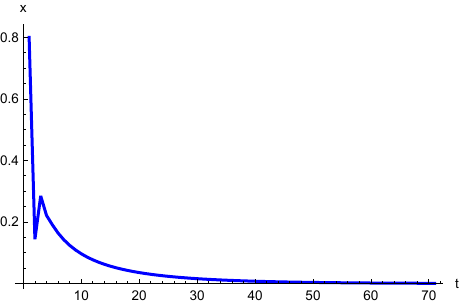}} 
	\subfloat[$r=0.9, x(0)=-0.2$]{\includegraphics[width = 2in]{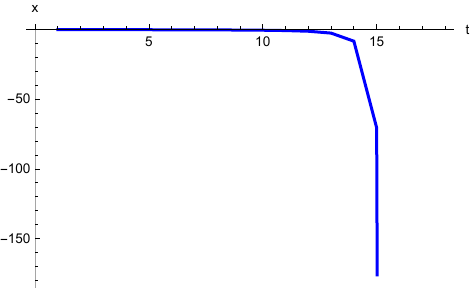}}\\
	\subfloat[$r=1.7, x(0)=0.8$]{\includegraphics[width = 2in]{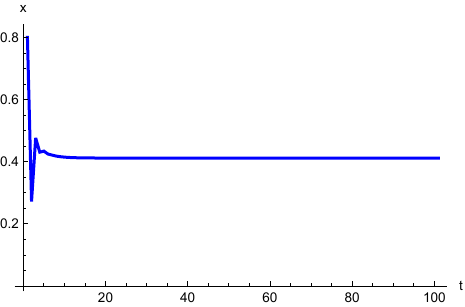}}
		\subfloat[$r=2.62, x(0)=0.8$]{\includegraphics[width = 2in]{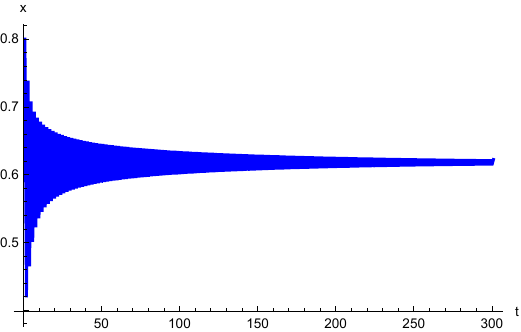}}\\
	\subfloat[$r=2.8, x(0)=0.8$]{\includegraphics[width = 2in]{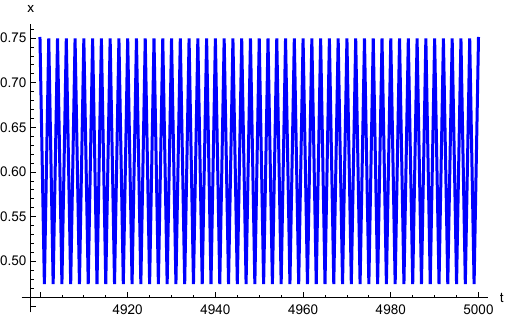}} 
	\subfloat[$r=3.3, x(0)=0.8$]{\includegraphics[width = 2in]{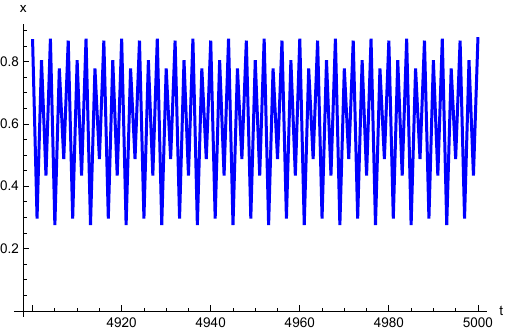}}\\
	\subfloat[$r=3.4, x(0)=0.8$]{\includegraphics[width = 2in]{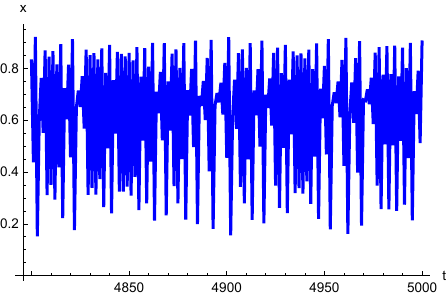}}
	\subfloat[$r=3.5, x(0)=0.8$]{\includegraphics[width = 2in]{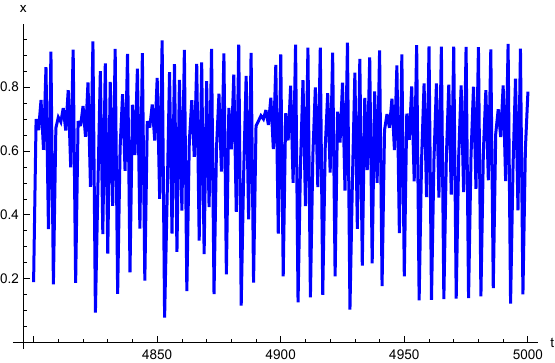}}\\
	\subfloat[$r=3.7, x(0)=0.8$]{\includegraphics[width = 2in]{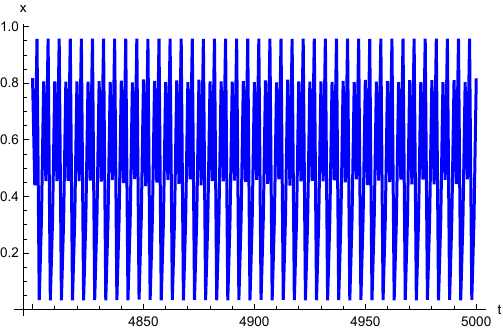}}
	\subfloat[$r=3.8, x(0)=0.8$]{\includegraphics[width = 2in]{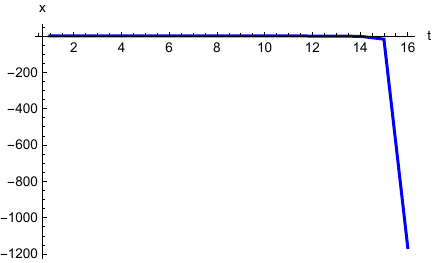}}	
	\caption{Solutions of the nonlinear system (\ref{nlin}) for $\alpha=0.6$, $q=0.3$ and various values of $r$}
	\label{fig6}
\end{figure}

For $r=0.9$, the equilibrium $x_{1*}=0$ is LAS and we get stable solutions for $x(0)$ close to $x_{1*}$ (cf. Fig. \ref{fig6} (a)) whereas unstable solutions near $x_{2*}=-0.1111$ (cf. Fig. \ref{fig6} (b)) as it is unstable.

For the fixed initial condition $x(0)=0.8$, we get stable solutions for $r\in (1, 2.62685)$ (see Fig. \ref{fig6} (c) with $r=1.7$, the solution is converging to the equilibrium $x_{2*}=0.411765$ and (d) with $r=2.62$). We observed asymptotically periodic solutions for $r\in [2.63, 3.3]$ (see Fig. \ref{fig6} (e) with $r=2.8$ and (f) with $r=3.3$ ). The chaotic oscillations are observed for $r\in [3.4, 3.6]$ (see Fig.  \ref{fig6} (g) with $r=3.4$ and (h) with $r=3.5$). The asymptotic periodic solutions are again observed for $r=3.7$ (Fig. \ref{fig6} (i)). The solutions becomes unbounded for $r=3.8$ (Fig. \ref{fig6} (j)).

A natural question arising from the previous construction is whether the proposed \(q\)-deformed fractional maps can be extended further by introducing multi-parameter deformations. In \(q\)-calculus and basic hypergeometric analysis, the \((p,q)\)-framework provides a broader generalization of classical \(q\)-analogues and has been successfully used in the study of special functions, combinatorics, and quantum algebras. Motivated by this, we now replace the Gaussian binomial coefficients appearing in the memory kernel by the more general \((p,q)\)-binomial coefficients. This extension introduces an additional deformation parameter, allowing greater flexibility in shaping the memory structure and stability geometry of the associated dynamical system. Moreover, the \(q\)-deformed fractional maps discussed in the previous sections are recovered as a special case when \(p=1\), showing that the present formulation naturally embeds the earlier theory into a more general deformation framework.

\section{\texorpdfstring{Generalization to include the $pq$-binomial coefficient}{Generalization to include the pq-binomial coefficient}}
	In this section, we generalize the map to include the $pq$-binomial coefficient.\\
The $pq$-binomial coefficient is defined as \cite{corcino2008p}:
\begin{equation}
	\binom{n}{k}_{pq}= \prod_{i=1}^{k} \frac{p^{n-i+1}-q^{n-i+1}}{p^i-q^{i}}.
\end{equation}
For $p=1$, this reduces to q-binomial coefficient.
\par In our case, we can define
\begin{equation}
	\tilde{\phi}_{\alpha,pq}(t) = \prod_{k=1}^{t} \frac{p^{\alpha
	+k-1}-q^{\alpha+k-1}}{p^{k}-q^{k}} = \binom{t+\alpha-1}{t}_{pq}.
\end{equation}
We have, 
\begin{align}
	\tilde{\phi}_{\alpha,pq}(t) &= \frac{\prod_{k=0}^{t-1} {p^\alpha p^{
		k}-q^\alpha q^{k}}}{\prod_{k=0}^{t-1}{p p^{k}-q q^{k}}} 
\nonumber \\
	&= \frac{((p^\alpha,q^\alpha);(p,q))_t}{((p,q);(p,q))_t}.
\end{align}

Furthermore, $(p,q)$-hypergeometric function is defined as \cite{jagannathan1998p}

\begin{eqnarray}
&&	_r\Phi_s((a_{1p},a_{1q});---; (p,q),z^{-1})\nonumber\\
	&&\quad =\sum_{n=0}^{\infty}\frac{((a_{1p},a_{1q});(p,q))_n [(-1)^n (q/p)^{n(n-1)/2}]^{1+s-r} z^{-n}}{((p,q);(p,q))_n},
\end{eqnarray}
where %$\vert q/p \vert <1$
 $r,s=0,1,2,\ldots$.\\

$\therefore$ For  $ r=1$ and $s=0$, we have
\begin{eqnarray}
&&	_1\Phi_0((p^\alpha,q^\alpha);---; (p,q),z^{-1})\nonumber\\
&&\quad	=\sum_{t=0}^{\infty}\frac{((p^\alpha,q^\alpha);(p,q))_t z^{-t}}{((p,q);(p,q))_t}\nonumber\\
&&\quad	=Z\left(\tilde{\phi}_{\alpha,pq}(t)\right).
\end{eqnarray}

%Consider the case $p=q^{-1}$ 
%\begin{equation}
%		\tilde{\phi}_{\alpha,q^{-1},q}(t) = 
%		\frac{(q^{2\alpha}; q^2)_t}
%	{(q^2;q^2)_t}
%	q^{t(1-\alpha)}.
%\end{equation}
%Therefore,
%\begin{eqnarray}
%Z\left(	\tilde{\phi}_{\alpha,q^{-1},q}(t)\right) &=& 
%\sum_{n=0}^{\infty}	\frac{(q^{2\alpha}; q^2)_n}{(q^2;q^2)_n}q^{n(1-\alpha)} z^{-n}\nonumber\\
%&=&\frac{(q^{2\alpha} q^{1-\alpha}z^{-1}; q^2)_\infty}
%{(q^{1-\alpha}z^{-1};q^2)_\infty}\nonumber\\
%&=& \frac{(q^{\alpha+1} z^{-1}; q^2)_\infty}
%{(q^{1-\alpha}z^{-1};q^2)_\infty},
%\end{eqnarray}
%where $\vert q^{1-\alpha}z^{-1}\vert <1$.

Consider a particular case $p=q^\rho$, where $\rho\in\Re$. We have
\begin{eqnarray}
	\tilde{\phi}_{\alpha,q^\rho,q}(t) &=& 
	\prod_{k=0}^{t-1} \frac{(q^{\rho(\alpha+k)}-q^{\alpha+k})}
	     {(q^{\rho(k+1)}-q^{k+1})}
	 = \prod_{k=0}^{t-1} 
	\frac{(1-q^{(\alpha+k)(1-\rho)})} {(1-q^{(k+1)(1-\rho)})}
	 q^{\rho(\alpha-1)} \nonumber\\
	 &=&  \frac{(q^{\alpha(1-\rho)}; q^{1-\rho})_t} {(q^{1-\rho};q^{1-\rho})_t}
	 q^{t\rho(\alpha-1)}.
\end{eqnarray}

%Here, $(a;q)_n$ denotes the $q$-Pochhammer symbol:
%\[ (a;q)_n = \prod_{k=0}^{n-1} (1 - aq^k) \].
Therefore,
\begin{eqnarray}
	Z\bigl(\tilde{\phi}_{\alpha,q^{\rho}q}(t)\bigr)
	&=& \sum_{n=0}^{\infty}\frac{(q^{\alpha(1-\rho)};q^{1-\rho})_n}
	{(q^{1-\rho};q^{1-\rho})_n}z^{-n} q^{n\rho(\alpha-1)} 
	= \frac{(q^{\alpha(1-\rho)}q^{\rho(\alpha-1)}z^{-1};q^{1-\rho})_\infty}
	{(z^{-1}q^{\rho(\alpha-1)};q^{1-\rho})_\infty} \nonumber\\
	&=& \frac{(q^{\alpha-\rho}z^{-1};q^{1-\rho})_\infty}
	{(q^{\rho(\alpha-1)}z^{-1};q^{1-\rho})_\infty}. \label{zt}
\end{eqnarray}

Using the $q$-binomial theorem
\[
\sum_{n=0}^{\infty}
\frac{(a;q)_n}{(q;q)_n}x^n
=
\frac{(ax;q)_\infty}{(x;q)_\infty},
\]
the representation (\ref{zt}) is valid provided
\[
|q^{1-\rho}|<1
\qquad\text{and}\qquad
|q^{\rho(\alpha-1)}z^{-1}|<1.
\]
Since $0<q<1$, the first condition is equivalent to
\[
\rho<1.
\]
For the stability analysis we consider $|z|=1$, and therefore the second condition reduces to
$
q^{\rho(\alpha-1)}<1.
$
Since $0<\alpha<1$, this is equivalent to
$
\rho<0.
$
Consequently, the stability analysis developed below is restricted to the regime
$
\rho<0,
$
or equivalently,
\[
p=q^\rho>1.
\]
We show in later section that for $p<1$, the memory kernel grows exponentially rather than decay, and is expected to diverge.

With these preliminaries, we are ready to consider the map

\begin{equation}
	x(t+1)=x(0)+(a-1) \bigl(\tilde{\phi}_{\alpha,pq} *x\bigr) (t), \label{genmap}
\end{equation}

	$0<\alpha<1, p=q^\rho, \vert q\vert <1 $ and $t \in \mathbb{N}_0$.

	Applying Z-transform we get,

	\begin{equation}
		zX(z)-zx(0)={\frac{x(0)}{1-z^{-1}} +(a-1) 
	 \frac{(q^{\alpha-\rho}z^{-1};q^{1-\rho})_\infty}
		{(q^{\rho(\alpha-1)}z^{-1};q^{1-\rho})_\infty}}X(z).
	\end{equation}

	$\therefore$ The characteristic equation becomes 

	\begin{equation}
		a =1+ \frac{(q^{\rho(\alpha-1)}z^{-1};q^{1-\rho})_\infty}{z(q^{\alpha-\rho}z^{-1};q^{1-\rho})_\infty}.
	\end{equation}

We obtain the boundary of stable region by substituting $z=e^{\iota \theta},$ $0\leq \theta \leq 2\pi$. i.e.
	\begin{equation}
		a =1+ \frac{(q^{\rho(\alpha-1)}e^{-\iota \theta};q^{1-\rho})_\infty}{e^{\iota \theta}(q^{\alpha-\rho}e^{-\iota \theta};q^{1-\rho})_\infty}. \label{stbdry}
	\end{equation}
    
	If $a\in \mathbb{R}$ then the bounds are given by $z=\pm 1$.

	The left and right bounds of the stable region on real line are given by
	\begin{equation}
		a_L =1- \frac{(-q^{\rho(\alpha-1)};q^{1-\rho})_\infty}{(-q^{\alpha-\rho};q^{1-\rho})_\infty} \label{lb}
	\end{equation}
	and 
	\begin{equation}
		a_R =1+ \frac{(q^{\rho(\alpha-1)};q^{1-\rho}))_\infty}{(q^{\alpha-\rho};q^{1-\rho})_\infty},\label{rb}
	\end{equation}
	respectively.

    For convenience, let us define
\begin{equation}
\Lambda_{\rho,q,\alpha}
=
\frac{(-q^{\rho(\alpha-1)};q^{\,1-\rho})_\infty}
     {(-q^{\alpha-\rho};q^{\,1-\rho})_\infty}; \Lambda'_{\rho,q,\alpha}=\frac{(q^{\rho(\alpha-1)};q^{1-\rho})_\infty}{(q^{\alpha-\rho};q^{1-\rho})_\infty}.
\end{equation}
Then the left and right stability bounds (\ref{lb}) can be written compactly as
\begin{equation}
a_L=1-\Lambda_{\rho,q,\alpha}, \;\;{\rm{and}} \;\; a_R=1+\Lambda'{\rho,q,\alpha}.
\end{equation}
This notation also makes the connection with the previously studied
$q$-deformed case transparent. Indeed, for $p=1$ we have $\rho=0$,
and therefore
\begin{equation}
\Lambda_{0,q,\alpha}
=
\frac{(-q^{0(\alpha-1)};q)_\infty}
     {(-q^{\alpha};q)_\infty}
=
\frac{(-1;q)_\infty}
     {(-q^{\alpha};q)_\infty}
=
\Lambda_{q,\alpha}.
\end{equation}
Consequently,
$
a_L
=
1-\Lambda_{0,q,\alpha}
=
1-\Lambda_{q,\alpha},
$
which is precisely the left stability bound obtained earlier for the
$q$-deformed fractional map. Moreover,
\begin{equation}
\Lambda'(\rho,q,\alpha)=\frac{(q^{\rho(\alpha-1)};q^{1-\rho})_\infty}
     {(q^{\alpha-\rho};q^{1-\rho})_\infty}
\Bigg|_{\rho=0}
=
\frac{(1;q)_\infty}
     {(q^\alpha;q)_\infty}
=0,
\end{equation}
since $(1;q)_\infty=0$. Hence the right stability bound reduces to
$
a_R=1.
$ Hence the $(p,q)$-deformed framework
contains the $q$-deformed theory as the special case $p=1$.

   Substituting $\rho=1$ (i.e. $p=q$) in the boundary (\ref{stbdry}) of stable region, we get
\begin{equation}
		a =1+e^{-\iota\theta}.
	\end{equation}
Thus, for $\rho=1$, the stable region is a circle centered at $(1,0)$ with
	radius 1.
\subsection{Example}
\begin{figure}
\centering
\includegraphics{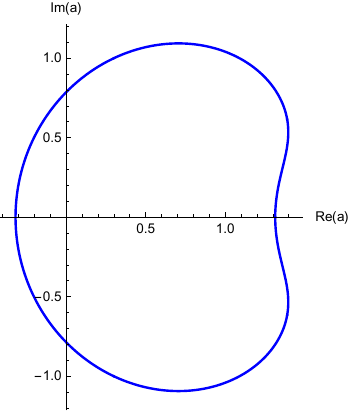}
\caption{Stable region for the system (\ref{genmap}) with $\alpha=0.4, q=0.8$ and $\rho=-1$}
\label{fig7}
\end{figure}

Consider the map defined by equation (\ref{genmap}) with $\alpha=0.4, q=0.8$ and $\rho=-1$ (i.e. $p=1/0.8$). Following the bounds (\ref{lb}) and (\ref{rb}), the system is asymptotically stable for $-0.318586<a<1.31415$. If the parameter $a$ is complex number then it should lie inside the boundary curve sketched in Figure \ref{fig7}.

\begin{figure}
	\subfloat[$a=-0.32$]{\includegraphics[width = 2in]{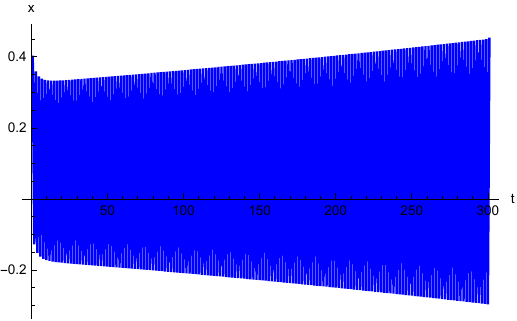}} 
	\subfloat[$a=1.32$]{\includegraphics[width = 2in]{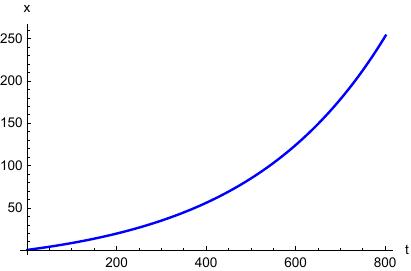}}\\
	\subfloat[$a=-0.31$]{\includegraphics[width = 2in]{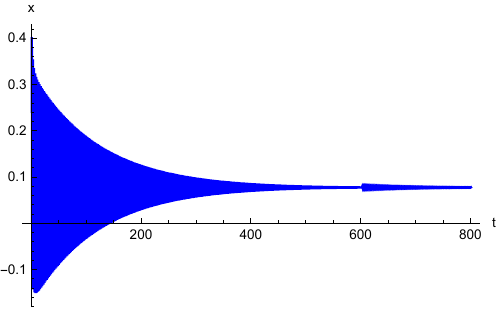}}
		\subfloat[$a=1.2$]{\includegraphics[width = 2in]{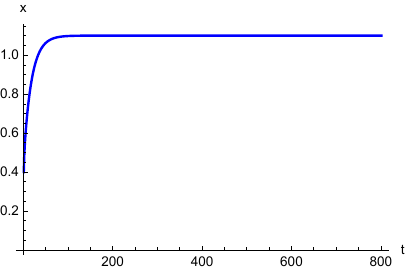}}\\
	\subfloat[$a=1.1+0.7\iota$]{\includegraphics[width = 2in]{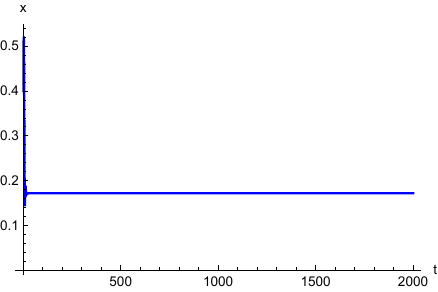}} 
	\subfloat[$a=1+1.1\iota$]{\includegraphics[width = 2in]{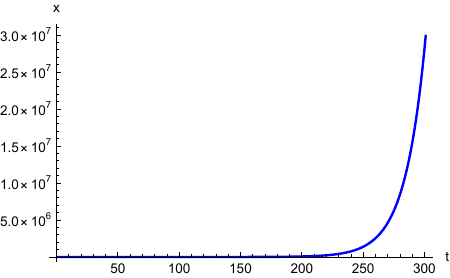}}	
	\caption{Solutions of the system (\ref{genmap}) for $\alpha=0.4, q=0.8$, $\rho=-1$ and various values of $a$}
	\label{fig8}
\end{figure}

For $a=-0.32$ and $a=1.32$ (outside the stable region), the system is unstable (cf. Figures \ref{fig8}(a) and (b)). For the parameter values $a=-0.31$ and $a=1.2$ (inside the stable region), the system is asymptotically stable (cf. Figures \ref{fig8}(c) and (d)). For the complex value $a=1.1+0.7\iota$ in the stable region, we get the convergent solution trajectory (Figure \ref{fig8}(e)). The unbounded solution for $a=1+1.1\iota$ outside stable region is sketched in Figure \ref{fig8}(f).

\begin{figure}
    \centering
    \includegraphics[width=5in]{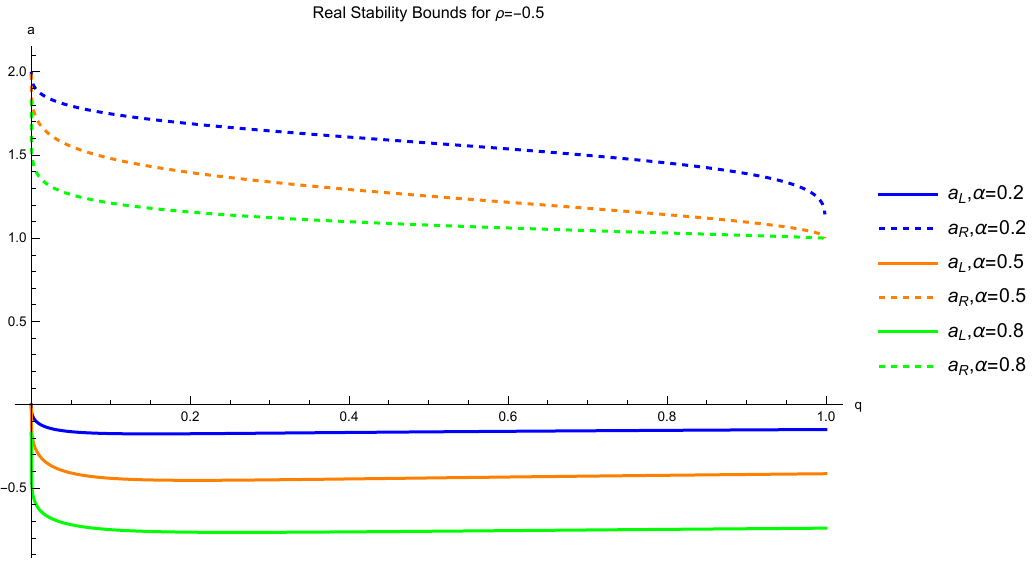}
    \caption{The real stability bounds $a_L$ (solid curves) and $a_R$ (dashed curves) for the $(p,q)$-deformed map with $\rho=-0.5$ ($p=q^{-0.5}$). Curves of the same color correspond to the same fractional order $\alpha$, with $\alpha=0.2$, $0.5$, and $0.8$. The admissible stability interval $a_L<a<a_R$ varies continuously with the deformation parameter $q$, while increasing $\alpha$ shifts both bounds upwards.
}
\label{fig:pqlim}
\end{figure}

Figure \ref{fig:pqlim} shows the dependence of the real stability bounds $a_L$ and $a_R$ on the deformation parameter $q$ for $\rho=-0.5$ and several values of $\alpha$. Both stability bounds vary continuously with $q$, while the location of the admissible interval $a_L<a<a_R$ is significantly influenced by the fractional order.

\subsection{Decay of the memory kernels}

An important feature of the proposed framework is the qualitative change
in the asymptotic behavior of the memory kernel as one moves from the
classical fractional case to the $q$- and $(p,q)$-deformed settings.

\paragraph{Classical fractional kernel: power-law decay.}

For the classical fractional order map, the kernel is
\[
\widetilde{\phi}_{\alpha}(n)
=
\binom{n+\alpha-1}{n}
=
\frac{\Gamma(n+\alpha)}
{\Gamma(\alpha)\Gamma(n+1)}.
\]

Using the asymptotic formula
\[
\frac{\Gamma(n+\alpha)}{\Gamma(n+1)}
\sim
n^{\alpha-1},
\qquad \text{as}  \,\,\, n\to\infty,
\]
we obtain
\[
\widetilde{\phi}_{\alpha}(n)
\sim
\frac{n^{\alpha-1}}{\Gamma(\alpha)}.
\]

Since $0<\alpha<1$, we have $\alpha-1\in(-1,0)$, and therefore
\[
\widetilde{\phi}_{\alpha}(n)
\propto
n^{-(1-\alpha)}.
\]

Thus, the classical fractional kernel exhibits algebraic power-law decay.
This slow decay is the characteristic signature of long-memory fractional
dynamics.

\paragraph{$q$-deformed kernel: asymptotic saturation.}

For the $q$-deformed case, the kernel is
$\widetilde{\phi}_{\alpha,q}(n)
=
\frac{(q^\alpha;q)_n}{(q;q)_n}=
\prod_{k=1}^{n}
\frac{1-q^{\alpha+k-1}}
{1-q^k}
$. 
We analyze its asymptotic behavior.

 Let,
$
x=q^k$. Hence
$$
\frac{1-q^{\alpha+k-1}}
{1-q^k}
=
\frac{1-q^{\alpha-1}x}
{1-x}
=
1+(1-q^{\alpha-1})q^k
+
O(q^{2k}).
$$

Using the expansion
$
\frac{1}{1-x}
=
1+x+O(x^2),
\; x\to 0.
$

Taking logarithms,
$$
\log \widetilde{\phi}_{\alpha,q}(n)
=
\sum_{k=1}^{n}
\log \left(
1+(1-q^{\alpha-1})q^k
+
O(q^{2k})
\right)
=
(1-q^{\alpha-1})
\sum_{k=1}^{n}q^k
+
O(1).
$$
Since
$
\log(1+y)=y+O(y^2)
$.

The geometric series
$
\sum_{k=1}^{\infty}q^k
=
\frac{q}{1-q}
$
converges for $0<q<1$, and hence the logarithmic series converges. 
Therefore,
$
\lim_{n\to\infty}
\widetilde{\phi}_{\alpha,q}(n)
=
\frac{(q^\alpha;q)_\infty}
{(q;q)_\infty}
=
C_q>0
$. 
More precisely,
$
\widetilde{\phi}_{\alpha,q}(n)
=C_q+O(q^n).
$
Thus the $q$-deformed kernel approaches a finite positive constant exponentially fast.

\paragraph{Remark.}
The convergence of the kernel to a nonzero constant may at first seem incompatible with stability, since the memory weights do not decay. However, stability is governed not only by the kernel but also by the feedback coefficient $(a-1)$ appearing in the map. Replacing the kernel asymptotically by its limiting value $C_q$, the memory term behaves approximately as
$
C_q(a-1)\sum_{j=0}^{n}x(j).
$
Introducing the cumulative sum
$
S_n=\sum_{j=0}^{n}x(j),
$
one obtains the asymptotic relation
$
S_{n+1}
\approx
\bigl(1+C_q(a-1)\bigr)S_n+x(0).
$
Hence the long-time behavior is controlled by the effective multiplier
$
\lambda_{\mathrm{eff}}=
1+C_q(a-1).
$
When $a<1$, the factor $(a-1)$ is negative and acts as a stabilizing feedback. Successive iterates may  alternate in sign, leading to substantial cancellations. Consequently, despite the kernel approaching a nonzero constant, the memory term need not grow without bound, and convergence  may still occur. For $q= 0$ case, we have $C_q= 1$ (See Eq. \ref{20}) and stability region is a circle. For finite $q$, it is  a constant only asymptotically and stability region is correspondingly deformed.  Still, an argument analogous to that leading to (\ref{22}) provides a heuristic explanation for the persistence of a nontrivial stability region. For $q\rightarrow 1$, $C_q\rightarrow 0$ and fractional behavior takes over leading to cardioid shaped stability region. Thus we may observe crescent shape for intermediate $q$ values.

Unlike the classical fractional kernel, which exhibits an algebraic power-law tail, the $q$-deformed kernel loses its scale-free character and instead approaches a finite asymptotic value with exponentially small corrections. Thus the $q$-deformation replaces long-range fractional memory by an asymptotically saturated memory profile.

\paragraph{$(p,q)$-deformed kernel: exponential decay, saturation, and growth.}

For the $(p,q)$-deformed case with $p=q^\rho$, the kernel is
\[
\widetilde{\phi}_{\alpha,q^\rho,q}(n)
=
\frac{(q^{\alpha(1-\rho)};q^{1-\rho})_n}
{(q^{1-\rho};q^{1-\rho})_n}
\,q^{n\rho(\alpha-1)}.
\]

Let
\[
r=q^{1-\rho}.
\]

Then
\[
\widetilde{\phi}_{\alpha,q^\rho,q}(n)
=
\frac{(r^\alpha;r)_n}{(r;r)_n}
\,q^{n\rho(\alpha-1)}.
\]

Using the identity
\[
(a;r)_n
=
\frac{(a;r)_\infty}{(ar^n;r)_\infty},
\]
we obtain
\[
\frac{(r^\alpha;r)_n}{(r;r)_n}
=
\frac{(r^\alpha;r)_\infty}{(r;r)_\infty}
\cdot
\frac{(r^{n+1};r)_\infty}
{(r^{n+\alpha};r)_\infty}.
\]

Since $0<r<1$,
$
r^{n+\alpha}\to0,
\qquad
r^{n+1}\to0,
\qquad \text{as} \quad  n\to\infty,
$
and therefore
$
(r^{n+\alpha};r)_\infty\to1,
\qquad
(r^{n+1};r)_\infty\to1.
$

Hence,
\[
\frac{(r^\alpha;r)_n}{(r;r)_n}
\to
\frac{(r^\alpha;r)_\infty}{(r;r)_\infty}
=
C_r.
\]

Thus the asymptotic behavior is governed by the factor
\[
q^{n\rho(\alpha-1)}
=
\exp\!\bigl(\rho(\alpha-1)n\log q\bigr).
\]

Since $0<q<1$, we have $\log q<0$, while $0<\alpha<1$ implies
$\alpha-1<0$.

Therefore:
\begin{itemize}

\item if $\rho<0$ (equivalently $p>1$), then
\[
q^{n\rho(\alpha-1)}
=p^{n(\alpha-1)}=
e^{-Dn},
\qquad D>0,
\]
and the kernel decays exponentially;

\item if $\rho=0$ (equivalently $p=1$), then the kernel approaches a
finite constant;

\item if $0<\rho<1$ (equivalently $q<p<1$), then the kernel grows
exponentially.

\end{itemize}

Consequently,
\[
\widetilde{\phi}_{\alpha,pq}(n)
\sim
C_r e^{-Dn}
\qquad (\rho<0),
\]
whereas for $0<\rho<1$ the kernel exhibits exponential growth. We have assumed $p\ge 1 $ in our case.

\begin{figure}
\includegraphics[width=5in]{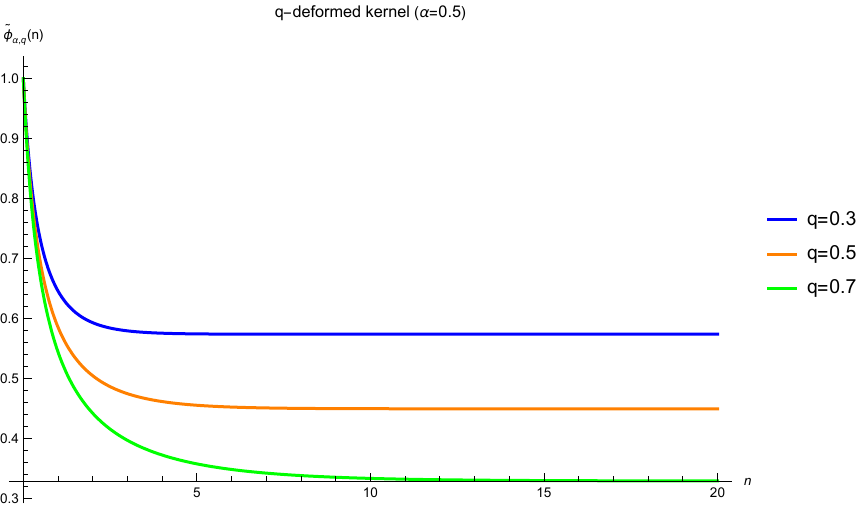}
\caption{The $q$-deformed kernel $\widetilde{\phi}_{\alpha,q}(n)$ for $\alpha=0.5$ and $q=0.3,0.5,0.7$. The kernels approach positive asymptotic constants $C_q$, confirming the asymptotic saturation predicted by the theory. }
\label{fig:qdef}
\end{figure}

\begin{figure}
\includegraphics[width=5in]{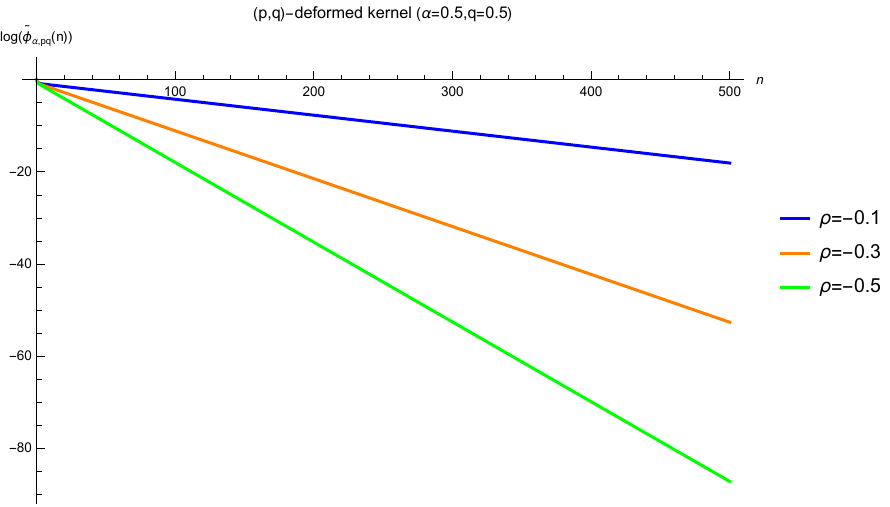}
\caption{The $(p,q)$-deformed kernel $\widetilde{\phi}_{\alpha,pq}(n)$ for $\alpha=0.5$, $q=0.5$, and $\rho=-0.1,-0.3,-0.5$ ($p=q^\rho>1$), plotted on a semi-logarithmic scale. The approximately linear behavior confirms the asymptotic relation $\widetilde{\phi}_{\alpha,pq}(n)\sim C_r e^{-Dn}$ and hence exponential decay of the memory kernel. More negative values of $\rho$ correspond to faster decay rates.}
\label{fig:pqdef}
\end{figure}

Thus the $(p,q)$-deformation introduces qualitatively different memory
regimes depending on the relative values of $p$ and $q$. The classical
power-law memory structure is recovered only in the singular limit
$q\to1^{-}$, where the $q$-Pochhammer symbols reduce to the classical
Gamma-function formulation.
Figure \ref{fig:qdef} illustrates the asymptotic saturation of the $q$-deformed memory kernel for several values of $q$, while Figure \ref{fig:pqdef} demonstrates the exponential decay of the $(p,q)$-deformed kernel when $\rho<0$ ($p>1$).

\section{Conclusion} \label{con.}

In this work, we have introduced a class of $q$-deformed fractional order maps by replacing the classical binomial kernel in discrete fractional dynamics with Gaussian (or $q$-) binomial coefficients. This construction provides a unified framework that connects classical fractional order maps, $q$-calculus, and discrete dynamical systems with memory.

A key feature of the proposed formulation is that the memory kernel is naturally deformed by the parameter $q$, leading to a smooth interpolation between different dynamical regimes. In particular, the classical fractional order map is recovered in the limit $q \to 1^-$, while the limit $q \to 0$ yields a memoryless discrete map. This establishes a clear bridge between short-term and long-term memory effects within a single parameterized structure.

We analyzed the stability of the zero solution using the Z-transform and derived a characteristic equation involving $q$-Pochhammer symbols. The stability boundary in the complex parameter space is described by a parametric curve whose geometry depends on both the fractional order $\alpha$ and the deformation parameter $q$. Numerical investigations demonstrate that this boundary undergoes a remarkable geometric transition from a circle (classical case) to a crescent-shaped curve for intermediate values of $q$, and approaches a cardioid in the fractional limit. This geometric deformation provides a clear visual interpretation of how $q$-deformation influences stability regions in fractional dynamical systems.

The framework was further extended to nonlinear maps, including a $q$-fractional logistic map, where linearization around equilibrium points yields stability conditions consistent with the linear theory. Numerical simulations confirm the analytical predictions and illustrate a rich variety of dynamical behaviors, including convergence to equilibria, periodic oscillations, and chaotic dynamics depending on system parameters.

We extend the proposed framework to \((p,q)\)-deformations through the introduction of \((p,q)\)-binomial coefficients in the memory kernel. This generalization incorporates an additional deformation parameter, providing greater flexibility in modeling discrete memory effects and stability structures. The resulting formulation contains the \(q\)-deformed fractional maps as a special case when \(p=1\), thereby embedding the present theory into a broader deformation framework. The appearance of \((p,q)\)-hypergeometric structures in the associated Z-transforms further indicates deep connections between fractional dynamics, special functions, and generalized combinatorial analysis. The present results suggest a broader hierarchy of memory kernels ranging from classical power-law memory to exponential memory induced by \((p,q)\)-deformations. It would be interesting to investigate whether scale-dependent deformations can generate stretched-exponential kernels and the corresponding stability geometries.

Overall, the present work establishes a broad class of deformed fractional discrete systems and reveals a deep interplay between memory effects, special functions, and geometric structures in stability analysis.
\section*{Declarations}
The authors declare that they have no known competing financial interests or personal relationships that could have appeared to influence the work reported in this paper.\\
Generative AI tools were used only to improve the language and 
clarity of this manuscript. The authors carefully reviewed and 
edited the output and take full responsibility for the final content.

	%\section*{Conflict of interest}
	
	%The authors declare that they have no conflict of interest.
	
	%	\end{itemize}


\bibliographystyle{aipnum4-2}      % AIP numerical bibliography style
\bibliography{ref.bib}

@book{gasper2004basic,
  title={Basic hypergeometric series},
  author={Gasper, George and Rahman, Mizan},
  volume={96},
  year={2004},
  publisher={Cambridge university press}
}

@misc{weisstein2017binomial,
  title={Binomial Coefficient, from MathWorld--A Wolfram Web Resource},
  author={Weisstein, Eric W},
  year={2017}
}

@article{luo2021fractional,
  title={Fractional chaotic maps with q--deformation},
  author={Luo, Cheng and Liu, Bao-Qing and Hou, Hu-Shuang},
  journal={Applied Mathematics and Computation},
  volume={393},
  pages={125759},
  year={2021},
  publisher={Elsevier}
}

@article{ran2022dynamics,
  title={On the dynamics of fractional q-deformation chaotic map},
  author={Ran, Jie and Li, Yu-Qin and Xiong, Yi-Bin},
  journal={Applied Mathematics and Computation},
  volume={424},
  pages={127053},
  year={2022},
  publisher={Elsevier}
}

@article{malik2025dynamical,
  title={Dynamical analysis, stabilization, and synchronization in the q-deformed discrete fractional Stefanski map.},
  author={Malik, MG and Awais, Muhammad and Ullah, Irfan and Bashir, Zia},
  journal={International Journal of Dynamics \& Control},
  volume={13},
  number={5},
  pages={1},
  year={2025}
}

@article{bhalekar2025dynamical,
  title={Dynamical analysis of fractional order generalized logistic map},
  author={Bhalekar, Sachin and Chevala, Janardhan and Gade, Prashant M},
  journal={Computational Mathematics and Mathematical Physics},
  volume={65},
  number={2},
  pages={424--441},
  year={2025},
  publisher={Springer}
}

@misc{mukhinsymmetric,
  title={Symmetric Polynomials and Partitions},
  author={Mukhin, Eugene},
  howpublished={Lecture notes, Indiana University--Purdue University Indianapolis},
  note={Unpublished lecture notes},
  year={2006}
}

@article{wu2020fractional,
  title={Fractional q-deformed chaotic maps: A weight function approach},
  author={Wu, Guang-Chao and Cankaya, M Niyazi and Banerjee, Santo},
  journal={Chaos: An Interdisciplinary Journal of Nonlinear Science},
  volume={30},
  number={12},
  pages={121106},
  year={2020},
  publisher={AIP Publishing}
}

@book{andrews1986q,
  title={q-Series: Their development and application in analysis, number theory, combinatorics, physics, and computer algebra},
  author={Andrews, George E},
  volume={66},
  year={1986},
  publisher={American Mathematical Soc.}
}

@book{hilfer2000applications,
  title={Applications of fractional calculus in physics},
  author={Hilfer, Rudolf},
  year={2000},
  publisher={World scientific}
}

@book{mainardi2022fractional,
  title={Fractional calculus and waves in linear viscoelasticity: an introduction to mathematical models},
  author={Mainardi, Francesco},
  year={2022},
  publisher={World Scientific}
}

@article{wu2014chaos,
  title={Chaos synchronization of the discrete fractional logistic map},
  author={Wu, Guo-Cheng and Baleanu, Dumitru},
  journal={Signal processing},
  volume={102},
  pages={96--99},
  year={2014},
  publisher={Elsevier}
}

@inproceedings{miller1989fractional,
  title={Fractional difference calculus},
  author={Miller, Kenneth S and Ross, Bertram},
  booktitle={Proceedings of the international symposium on univalent functions, fractional calculus and their applications},
  pages={139--152},
  year={1989}
}

@article{atici2009initial,
  title={Initial value problems in discrete fractional calculus},
  author={Atici, Ferhan and Eloe, Paul},
  journal={Proceedings of the American Mathematical Society},
  volume={137},
  number={3},
  pages={981--989},
  year={2009}
}

@article{mozyrska2015transform,
  title={The-transform method and delta type fractional difference operators},
  author={Mozyrska, Dorota and Wyrwas, Ma{\l}gorzata and others},
  journal={Discrete Dynamics in Nature and Society},
  volume={2015},
  year={2015},
  publisher={Hindawi}
}

@book{ernst2000history,
  title={The History of Q-calculus and a New Method},
  author={Ernst, T.},
  series={UUDM Report},
  url={https://books.google.co.in/books?id=pnIntwAACAAJ},
  year={2000},
  publisher={Department of Mathematics, Uppsala University}
}

@book{koepf1998hypergeometric,
  title={Hypergeometric summation},
  author={Koepf, Wolfram},
  year={1998},
  publisher={Springer}
}

@article{konvalina2000unified,
  title={A unified interpretation of the binomial coefficients, the Stirling numbers, and the Gaussian coefficients},
  author={Konvalina, John},
  journal={The American Mathematical Monthly},
  volume={107},
  number={10},
  pages={901--910},
  year={2000},
  publisher={Taylor \& Francis}
}

@article{yamano2002some,
  title={Some properties of q-logarithm and q-exponential functions in Tsallis statistics},
  author={Yamano, Takuya},
  journal={Physica A: Statistical Mechanics and its Applications},
  volume={305},
  number={3-4},
  pages={486--496},
  year={2002},
  publisher={Elsevier}
}

@book{podlubny1998fractional,
  title={Fractional differential equations: an introduction to fractional derivatives, fractional differential equations, to methods of their solution and some of their applications},
  author={Podlubny, Igor},
  year={1998},
  publisher={Elsevier}
}

@book{kilbas2006theory,
  title={Theory and applications of fractional differential equations},
  author={Kilbas, Anatoli{\u\i} Aleksandrovich and Srivastava, Hari M and Trujillo, Juan J},
  volume={204},
  year={2006},
  publisher={elsevier}
}

@article{atici2007transform,
  title={A transform method in discrete fractional calculus},
  author={Atici, Ferhan M and Eloe, Paul W},
  journal={International Journal of Difference Equations},
  volume={2},
  number={2},
  year={2007}
}

@article{pakhare2022synchronization,
  title={Synchronization in coupled integer and fractional-order maps},
  author={Pakhare, Sumit S and Bhalekar, Sachin and Gade, Prashant M},
  journal={Chaos, Solitons \& Fractals},
  volume={156},
  pages={111795},
  year={2022},
  publisher={Elsevier}
}

@article{joshi2023controlling,
  title={Controlling fractional difference equations using feedback},
  author={Joshi, Divya D and Bhalekar, Sachin and Gade, Prashant M},
  journal={Chaos, Solitons \& Fractals},
  volume={170},
  pages={113401},
  year={2023},
  publisher={Elsevier}
}

@article{bhalekar2023fractional,
  title={Fractional-order periodic maps: Stability analysis and application to the periodic-2 limit cycles in the nonlinear systems},
  author={Bhalekar, Sachin and Gade, Prashant M},
  journal={Journal of Nonlinear Science},
  volume={33},
  number={6},
  pages={119},
  year={2023},
  publisher={Springer}
}

@article{gade2021fractional,
  title={On fractional order maps and their synchronization},
  author={Gade, Prashant M and Bhalekar, Sachin},
  journal={Fractals},
  volume={29},
  number={06},
  pages={2150150},
  year={2021},
  publisher={World Scientific}
}

@article{edelman2023stability,
	title={Stability of fixed points in generalized fractional maps of the orders 0< $\alpha$< 1},
	author={Edelman, Mark},
	journal={Nonlinear Dynamics},
	volume={111},
	number={11},
	pages={10247--10254},
	year={2023},
	publisher={Springer}
}

@article{bhalekar2022stability,
	title={Stability analysis of fixed point of fractional-order coupled map lattices},
	author={Bhalekar, Sachin and Gade, Prashant M},
	journal={Communications in Nonlinear Science and Numerical Simulation},
	volume={113},
	pages={106587},
	year={2022},
	publisher={Elsevier}
}

@article{samko1993integration,
	title={Integration and differentiation to a variable fractional order},
	author={Samko, Stefan G and Ross, Bertram},
	journal={Integral transforms and special functions},
	volume={1},
	number={4},
	pages={277--300},
	year={1993},
	publisher={Taylor \& Francis}
}

@article{bhalekar2022stabilityc,
	title={Stability and dynamics of complex order fractional difference equations},
	author={Bhalekar, Sachin and Gade, Prashant M and Joshi, Divya},
	journal={Chaos, Solitons \& Fractals},
	volume={158},
	pages={112063},
	year={2022},
	publisher={Elsevier}
}

@article{corcino2008p,
  title={On p, q-binomial coefficients},
  author={Corcino, Roberto B and On, P},
  journal={Integers},
  volume={8},
  number={1},
  pages={29},
  year={2008}
}

@article{jagannathan1998p,
  title={(P, Q)-special functions},
  author={Jagannathan, R},
  journal={arXiv preprint math/9803142},
  year={1998}
}

@article{jackson1909q,
  title={On q-functions and a certain difference operator},
  author={Jackson, Frank Hilton},
  journal={Transactions of the Royal Society of Edinburgh},
  volume={46},
  pages={253--281},
  year={1909}
}

@book{kac2002quantum,
  title={Quantum Calculus},
  author={Kac, Victor and Cheung, Pokman},
  year={2002},
  publisher={Springer}
}

@book{annaby2012q,
  title={q-Fractional Calculus and Equations},
  author={Annaby, Mahmoud H. and Mansour, Zeinab S.},
  year={2012},
  publisher={Springer}
}

@book{goodrich2015discrete,
  title={Discrete Fractional Calculus},
  author={Goodrich, Christopher S. and Peterson, Allan C.},
  year={2015},
  publisher={Springer}
}

\end{document}